\documentclass[aop]{imsart}

\RequirePackage{amsthm,amsmath,amsfonts,amssymb}
\RequirePackage[numbers]{natbib}
\RequirePackage[colorlinks,citecolor=blue,urlcolor=blue]{hyperref}

\startlocaldefs
\numberwithin{equation}{section}
\theoremstyle{plain}
\newtheorem{theo}{Theorem}[section]
\newtheorem{prop}{Proposition}[section]
\newtheorem{lem}[prop]{Lemma}

\theoremstyle{remark}
\newtheorem{remark}[prop]{Remark}

\newcommand{\N}{\mathbb{N}}
\newcommand{\Z}{\mathbb{Z}}
\newcommand{\Q}{\mathrm{Q}} 
\newcommand{\R}{\mathbb{R}}

\def \1{\mathbf{1}}

\def\({\left(}
\def\){\right)}

\def\nab{\nabla}
\def\pa{\partial}
\def \ol#1{\overline{#1}}
\def\inner#1{\left \langle #1 \right \rangle}

\def\car{\Box}
\def \Var{\operatorname{Var}}
\def \Cov{\operatorname{Cov}}
\def \supp{\operatorname{supp}} 
\newcommand{\Esp}{\mathbb{E}}
\def \dist{\mathrm{dist}}

\def \MM{\mathrm{MM}}
\def\cl{\mathrm{cl}}
\def\K{\mathsf{K}}
\def \epsilon{\varepsilon}

\newcommand{\tref}[1]{Theorem~\ref{t.#1}}
\newcommand{\pref}[1]{Proposition~\ref{p.#1}}
\newcommand{\lref}[1]{Lemma~\ref{l.#1}}
\newcommand{\cref}[1]{Corollary~\ref{c.#1}}

\newcommand{\sref}[1]{Section~\ref{s.#1}}
\newcommand{\rref}[1]{Remark~\ref{r.#1}}
\newcommand{\eref}[1]{(\ref{e.#1})}

\endlocaldefs

\begin{document}
	
	\begin{frontmatter}
		\title{Thermodynamic and Scaling Limits of the non-Gaussian Membrane Model }
		\runtitle{Limits of the non-Gaussian Membrane Model}
		
		\begin{aug}
			\author[A]{\fnms{Eric} \snm{Thoma}\ead[label=e1]{eric.thoma@cims.nyu.edu}},
			
			\address[A]{Courant Institute, New York University, \printead{e1}}
		\end{aug}
		
		\begin{abstract}
			We characterize the behavior of a random discrete interface $\phi$ on $[-L,L]^d \cap \Z^d$ with energy $\sum V(\Delta \phi(x))$ as $L \to \infty$, where $\Delta$ is the discrete Laplacian and $V$ is a uniformly convex, symmetric, and smooth potential. The interface $\phi$ is called the non-Gaussian membrane model. By analyzing the Helffer-Sj\"ostrand representation associated to $\Delta \phi$, we provide a unified approach to continuous scaling limits of the rescaled and interpolated interface in dimensions $d=2,3$, Gaussian approximation in negative regularity spaces for all $d \geq 2$, and the infinite volume limit in $d \geq 5$. Our results generalize some of those of \cite{CDH19}.
		\end{abstract}
		
		
		
	\end{frontmatter}
	
	\begin{section}{Introduction} \label{s.intro_sec}
		In the present article we investigate the membrane model, which is a random scalar field $(\phi(x))_{x \in \Z^d}$ with distribution given by
		\begin{equation*}
			\MM_L(d\phi) = \frac{1}{\mathrm{Z}(\car_L)} \exp\(-\sum_{x \in \Z^d} V(\Delta \phi(x)) \) \prod_{x \in \car_L} d\phi(x) \prod_{x \in \Z^d \setminus \car_L} \delta_0(d\phi(x)),
		\end{equation*}
		where $\car_L = [-L,L]^d \cap \Z^d$ is a $d$-dimensional discrete cube, $V : \R \to \R$ is a potential satisfying certain conditions, and $\Delta$ is the discrete Laplacian on $\Z^d$. The constant $\mathrm{Z}(\car_L)$ is a normalizing factor, and we have enforced $\phi(x) = 0$ for $x \not \in \car_L$ through the Dirac delta measure $\delta_0$, though other boundary conditions and other domains besides cubes are possible.
		
		In particular, we will investigate limits of $\MM_L$ as $L \to \infty$ in $d \geq 2$. There are multiple senses in which we can take a limit, and we will consider three types. The first type is the infinite-volume limit, also known as the thermodynamic limit, which is the distributional limit of $(\phi(x))_{x \in \car_K}$ as $L \to \infty$ and $K$ stays fixed. This limit only exists in the supercritical dimensions $d \geq 5$, and is generally non-Gaussian. The second and third types concern the scaling limit of the rescaled field $\ol{\phi}(x) = L^{d/2 - 2} \phi(Lx)$ defined on $[-1,1]^d \cap \frac{1}{L} \Z^d$. We consider for $d \geq 2$ the random quantity 
		$$
		L^{-d} \sum_{x \in [-1,1]^d \cap \frac{1}{L} \Z^d} \ol{\phi}(x) f(x),
		$$
		where $f : [-1,1]^d \to \R$ is a smooth enough function, and approximate its law by a Gaussian. In dimensions $d=2,3$, we consider an appropriate interpolation of $\ol{\phi}(x)$ to a function with domain $[-1,1]^d$ and prove convergence of its distribution in the space of continuous functions.
		
		We prove the limits under the condition that $V$ is symmetric, uniformly convex, and $C^3$ with uniformly bounded second and third derivatives. Previously, similar limits were considered in \cite{CDH19} in the case where $V$ is quadratic, and therefore $\mathrm{MM}_L$ is Gaussian. Unfortunately, the proof methods cannot be easily generalized to the case of more general $V$. We introduce a new approach for the membrane model that allows us to treat all three types of limits in a non-Gaussian setting and under a unified framework. Our methods can also be adapted for regimes of general {\it semi-flexible membranes} in which the Laplacian energy dominates. In the case of limits of the rescaled field $\ol{\phi}$, we prove quantitative Gaussian approximation results and use the results of \cite{CDH19} to characterize the limit.
		
		We study the membrane model entirely through the distribution of the discrete Laplacian of the field $\phi$. We consider the equation
		\begin{equation*}\begin{cases}
				\Delta \phi(x) = \eta(x), \quad &x \in \car_L,\\
				\phi(x) = 0, \quad &x \in \partial \car_L,
		\end{cases} \end{equation*}
		where $\pa \car_L$ is the outer boundary of $\car_L$ in $\Z^d$, which gives a correspondence between $\phi$ and a field $\eta \in \R^{\car_L}$. Given the membrane model $\phi$ on $\car_L$, we can find $\eta$ by taking the discrete Laplacian, and given $\eta$, we can recover $\phi$ by solving the above Dirichlet problem. We let $\Q_L$ be the distribution of $\eta$ when $\phi$ is distributed by $\MM_L$. The measure $\Q_L$ is itself a Gibbs measure on $\R^{\car_L}$ with energy
		$$
		\sum_{x \in \car_L} V(\eta(x)) + \sum_{z \in \pa \car_L} V\(-\sum_{x \in \car_L} P^z_L(x) \eta(x) \)
		$$
		where $P^z_L(x)$ is the probability a simple random walk on $\Z^d$ started at $x$ exits $\car_L$ at $z \in \pa \car_L$, i.e.\ the Poisson kernel. Note that the latter sum in the energy is present due to the fact that $\Delta \phi$ is supported on $\car_L \cup \pa \car_L$, whereas, by definition, $\eta$ consists only of the data of $\Delta \phi$ on $\car_L$, and so we must recover the values of $\Delta \phi(z)$ for $z \in \pa \car_L$ from $\eta$ by using the constraint that $\phi$ is supported on $\car_L$. We use variational methods to understand the Helffer-Sj\"ostrand representation associated with $\eta$, which turns out to be more tractable than that associated with $\phi$. To our knowledge, the approach of studying $\Delta \phi$, and ultimately $\phi$, through the Gibbs measure $\Q_L$ is new, and we expect this approach may be fruitful for further investigation of the non-Gaussian membrane model.

		\begin{subsection}{Background and Motivation}
			The membrane model is an instance of a discrete interface model. These are random fields $\phi \in \R^{\Z^d}$ with law proportional to
			$$
			\exp(-H(\phi)) \prod_{x \in \car_L} d\phi(x) \prod_{x \in \Z^d \setminus \car_L} \delta_0(d\phi(x)),
			$$
			for a Hamiltonian $H : \R^{\Z^d} \to \R$. Perhaps the most studied of such models is the discrete Gaussian Free Field (GFF), which corresponds to the model with $H(\phi) \propto \sum_{x \in \Z^d} |\nab \phi(x)|^2$ where $\nab$ is the discrete gradient. Here, as with the membrane model, it is a natural question to consider the possible scalings of $\phi$ for which the $L \to \infty$ limit exists, and for the free field there are many powerful tools available arising from the Gaussian nature and the gradient form of the energy. Perhaps most notably, there is a random walk representation of the covariance of the GFF.
			
			More generally, a major research direction has been to extend results for the discrete Gaussian free field to the "$\nabla \phi$" model, which corresponds to $H(\phi) = \sum_{x \in \Z^d} V(\nab \phi(x))$ and $V$ uniformly convex, symmetric, and smooth. We refer to \cite{F05} for a survey of results for this model.
			
			Notably, in \cite{AW19} methods from elliptic stochastic homogenization were used to analyze the model (see also \cite{AD22}). The gradient form of the energy allows for methods from elliptic regularity, such as the De Giorgi-Nash-Moser theory, to be applied. We do not know of a replacement for these methods for the membrane model, and this is a major obstacle in applying a similar approach. The general method of \cite{AW19}, i.e.\ analyzing the Helffer-Sj\"ostrand equation, was first applied to the $\nab \phi$ model in \cite{NS97}, though the method's history goes back further (\cite{HS94}, \cite{S96}). In \cite{GOS01}, The Helffer-Sj\"ostrand equation for the $\nab \phi$ model is analyzed through a random walk in a dynamic random environment, a connection also unavailable in the present context. We will follow in using Helffer-Sj\"ostrand equations, but must use quite different methods of analysis.
			
			The physical interest in the membrane model stems partly from its role among the family of semi-flexible membranes, which are discrete interface models with
			$$
			H(\phi) = \sum_{x \in \Z^d}  V_1(\nab \phi(x)) + V_2(\Delta \phi(x))
			$$
			where $V_1$ and $V_2$ are potentials characterizing the membrane's lateral tension and bending rigidity (see \cite{CDH21} and references therein). The membrane model is also interesting due to its scaling properties, especially in the critical dimension $d=4$, which plays a role analogous to that dimension $2$ plays for the discrete GFF. In this dimension, correlations are logarithmic.
			
			We mention some relevant results about the Gaussian membrane model proved in \cite{CDH19}. The covariance operator of the Gaussian membrane model is the inverse discrete bi-Laplacian operator with "pinned" Dirichlet boundary conditions. An analysis of the scaling limit of this operator in $d=2,3$, proving convergence to the solution of a continuum problem, was achieved in \cite{MS19} and is used for the $d=2,3$ scaling limit results of \cite{CDH19}. Specifically, \cite{CDH19} proves that the rescaled and interpolated membrane model $\ol{\phi}$ converges weakly as $L \to \infty$ to a continuum Gaussian process on $[-1,1]^d$. In $d \geq 4$, \cite{CDH19} proves limits of $\ol{\phi}$ (on smooth domains instead of $\car_L$) using a framework of Gaussian fields and an extension of a method from \cite{T64} on the approximation of discrete equations by PDE.
			
			Besides for scaling limits, other questions of interest for the membrane model (and many other interface models) include entropic repulsion, pinning, wetting, behavior of the interface maximum, and level set percolation. Entropic repulsion was addressed in $d \geq 5$ by \cite{K07} and \cite{S03}, in $d=4$ by \cite{K09} (and the thesis \cite{K08}), and in $d=2,3$ by \cite{BDKS19}. For pinning in $d \geq 4$, results are given in \cite{S20a}, and pinning in $d =2,3$ is not well understood. The behavior of the maximum height of the membrane for the critical dimension $d=4$ was addressed in \cite{S20b} and for $d \geq 5$ in \cite{CCH16}. Level-set percolation was recently investigated in \cite{CN22} and \cite{M22}. These results are all for the Gaussian model. For the non-Gaussian model, there are very few results. Notably, \cite{CD09} proved the scaling limit for the non-Gaussian model in $d=1$ using a random walk representation specific to one-dimension.
		\end{subsection}
		
		\begin{subsection}{Statement of Main Results}
			We now state the main results of the paper. The parameter $L$ will be assumed to be a positive integer throughout. We assume that the potential $V$ is $C^3$ with $\sup V'' < \infty$, $\sup |V^{(3)}| < \infty$, symmetric about $0$, and $\inf V'' > 0$. \rref{potential.assumption} comments on possible extensions of our results to more general potentials.
			
			For the infinite volume limit, we characterize the non-Gaussian infinite volume limit of $\phi$ in $d \geq 5$. In $d \leq 4$, the model experiences a well-known blow up of variances as $L \to \infty$, and so the infinite volume limit does not exist.
			
			\begin{theo}[Infinite Volume Limit] \label{t.inf_limit_thm}
				Let $d \geq 5$, and let $-\Gamma : \Z^d \times \Z^d \to \R$ be the Green's function for the discrete Laplacian on $\Z^d$. Let $a \in \R^{\Z^d}$ have compact support. Let $\nu^0$ be a probability measure on $\R$ with Lebesgue density proportional to $\exp(-V(\xi))$, and let $\xi$ be a $\nu^0$ distributed random variable.
				
				Then we have
				\begin{equation} \label{e.inf_limit}
					\lim_{L \to \infty} \log \Esp_{\MM_L}  \exp\(\sum_{x \in \Z^d} a(x) \phi(x)\)  = \sum_{x \in \Z^d} \log \Esp_{\nu^0} \exp \( (\Gamma \ast a)(x)  \xi \)
				\end{equation}
				where $(\Gamma \ast a)(x) = \sum_{y \in \Z^d} \Gamma(x,y) a(y)$. In particular, for any fixed integer $K$, the field $(\phi(x))_{x \in \car_K}$ converges weakly and in moments to a random variable with cumulant generating function given by the RHS of \eref{inf_limit}.
			\end{theo}
			
			In the case that $V$ is quadratic, the RHS of \eref{inf_limit} is $\sum_{x \in \Z^d} \frac{1}{2 V''} |\Gamma \ast a(x)|^2.$ It follows that the limit is a Gaussian field with covariance operator $\frac12 (V'')^{-1} \Gamma \ast \Gamma$. Such a result is well-known; see e.g.\ \cite{K08}. For general $V$, the limiting distribution is non-Gaussian, since the RHS of \eref{inf_limit} is non-quadratic in $a$. Our result implies that the Laplacian field $\Delta \phi$ converges to an i.i.d. field, and the limiting distribution of each $\Delta \phi(x)$ is given by $\nu^0 \propto \exp(-V)$. The infinite volume law of $\phi$ is in fact the pushforward of $(\nu^0)^{\otimes \Z^d}$ under the convolution by $\Gamma$, which is well-defined due to the fact that $\Gamma(x,\cdot) \in L^2(\Z^d)$ in $d \geq 5$.
			
			Our second and third results involve the limit of the rescaled field $\ol{\phi}(x) = L^{d/2 - 2} \phi(L x)$, defined on a lattice of spacing $1/L$. In $d \geq 4$, the limit as $L \to \infty$ of $\ol{\phi}(x)$ should exist in a certain negative regularity sense, i.e.\ if we consider $\ol{\phi}$ summed against smooth test functions. In contrast to the infinite volume limit of \tref{inf_limit_thm}, the limit will be Gaussian due to the averaging of fluctuations in the summation. In $d=2,3$, the limit will be a H\"older continuous Gaussian process after an appropriate interpolation of the interface. The Gaussianity in $d=2,3$ comes from $\Gamma \not \in L^2(\Z^d)$, and so the field $\ol{\phi}$ is a large scale average of $\Delta \ol{\phi}$.
			
			In the case that $V$ is quadratic and $d=2,3$, the scaling limit was proved on squares and cubes in \cite{CDH19}. In $d \geq 4$, the scaling limit was proved on bounded smooth domains, i.e.\ with the family $\car_L$ replaced by $L \Omega \cap \Z^d$ for a bounded smooth domain $\Omega \subset \R^d$ and $L \Omega = \{Lx : x \in \Omega\}$.
			
			We need a technical assumption, which will allow us to prove in \sref{append} that certain fields do not concentrate on a boundary layer of $\car_L$. For a given $f : [-1,1]^d \to \R$ and $\Delta_{\R^d}$ the continuum Laplacian, consider the solution $u$ to
			\begin{equation} \label{e.regularity.assumption}
				\begin{cases}
					\Delta_{\R^d}^2 u(x) = f(x), \quad &x \in (-1,1)^d,\\
					u(x) = \pa_n u(x) = 0, \quad &x \in \pa (-1,1)^d.
				\end{cases}
			\end{equation}
			Our assumption is that $u \in C^5([-1,1]^d)$ (see \rref{regularity.assumption}).
			
			Our result approximates the cumulant generating function of the non-Gaussian model by a Gaussian model with an effective covariance and with explicit error estimates. We are in particular interested in the statistic
			$$
			\frac{1}{(2L)^d} \sum_{x \in L^{-1} \Z^d \cap [-1,1]^d} f(x) \ol{\phi}(x),
			$$
			which is a discrete approximation of the continuum $L^2([-1,1])$ inner product. We make $\ol{\phi}$ into an operator on smooth functions in this way, and we will approximate the law of $\ol{\phi}$ by a Gaussian law.
			
			\begin{theo}[Distributional Gaussian Approximation] \label{t.scal_limit_int}
				Let $d \geq 2$ and $f :[-1,1]^d \to \R$ be such that $u \in C^5([-1,1]^d)$ as defined above. Let $\ol{\phi}(x) = L^{d/2 - 2} \phi(L  x)$ be the rescaled membrane, defined on $[-1,1]^d \cap L^{-1} \Z^d$, and let $\MM_L^G$ be the (Gaussian) membrane model on $\car_L$ with potential $V^G(\xi) = \frac12 |\xi|^2$. Then we have
				\begin{align} \label{e.scal_limit_int}
					\lefteqn{\log \Esp_{\MM_L} \exp \( L^{-d} \sum_{x \in L^{-1} \Z^d \cap [-1,1]^d} f(x) \ol{\phi}(x)\)} \quad & \\ \notag &= \log \Esp_{\MM_L^G}  \exp \( (\Var_{\nu^0} \xi)^{1/2} L^{-d} \sum_{x \in L^{-1} \Z^d \cap [-1,1]^d} f(x) \ol{\phi}(x)\) + \mathrm{Error},
				\end{align}
				and $|\mathrm{Error}| \leq CL^{-\frac{d-1}{6d-2}}(1 + (\log L)^3\1_{d=2})$. The constant $C$ depends on $d$ and is a polynomial in $\| u \|_{C^5([-1,1]^d)}$, $(\inf V'')^{-1}$, $\sup V''$, and $\sup |V^{(3)}|$. The symbols $\nu_0$ and $\xi$ are the same as in \tref{inf_limit_thm}.
			\end{theo}
			
			While we have stated the above result for the membrane model on $\car_L$ and test functions on $[-1,1]^d$, we are also interested in the membrane model on $L\Omega \cap \Z^d$ and test functions on $\Omega$, for bounded, smooth domains $\Omega \subset \R^d$, particularly due to the importance of this setting in \cite{CDH19}. To that end, we mention that \tref{scal_limit_int} also holds, with minimal changes to the proof, when adapted to this setting (and even weaker conditions on $\Omega$ would suffice). In the smooth boundary case, the limit as $L \to \infty$ of the RHS of \eref{scal_limit_int} for sufficiently smooth $f$ compactly supported in $\Omega$ is computed in (\cite{CDH19}, Section 3).
			
			For the next result, we will need an interpolation of the rescaled membrane $\ol{\phi}$, defined initially on $[-1,1]^d \cap L^{-1} \Z^d$, to a continuous function on $[-1,1]^d$, for $d=2,3$. By an abuse of notation, we will also denote the interpolated interface by $\ol{\phi}$. We refer to (\cite{CDH19}, Section 2.1) for the specific interpolation procedure, though our main estimate \eref{scal_limit_pt_approx} holds under essentially any reasonable interpolation procedure.
			
			\begin{theo}[Scaling Limit in $d=2,3$]\label{t.scal_limit_pt}
				Let $d=2,3$ and let $\ol{\phi}(x) = L^{d/2 - 2} \phi(Lx)$ be the rescaled membrane on $L^{-1} \Z^d \cap [-1,1]^d$, and interpolate the interface to $[-1,1]^d$ as done in \cite{CDH19}. Let $x_1, x_2, \ldots, x_k$ be points in $[-1,1]^d$. Then for any $a \in \R^k$, we have
				\begin{equation} \label{e.scal_limit_pt_approx}
					\log \Esp_{\MM_L} \exp\(\sum_{i = 1}^k a_i \ol{\phi}(x_i)\) = \log \Esp_{\MM_L^G}  \exp \((\Var_{\nu^0} \xi)^{1/2} \sum_{i = 1}^k a_i \ol{\phi}(x_i)\) + \mathrm{Error}
				\end{equation}
				for $|\mathrm{Error}| \leq CL^{-\frac{d-1}{6d-2}}(1+ (\log L)^3\1_{d=2})$. Moreover, there exists a Gaussian process $\Psi : [-1,1]^d \to \R$ such that $\ol{\phi}$ converges in distribution to $\Psi$ in the space of continuous functions on $[-1,1]^d$ as $L \to \infty$. The process $\Psi$ is centered and a.s.\ $\alpha$-H\"older continous for every parameter $\alpha \in (0,1)$ in $d=2$ or $\alpha \in (0,1/2)$ in $d = 3$. 
			\end{theo}
			
			The above theorem not only concerns the limit of the finite dimensional distributions of $\ol{\phi}$, but also the distribution of the maximum of $\ol{\phi}$. The covariance of the limiting process $\Psi$ is proportional to the Green's function of the continuum bi-Laplacian on $[-1,1]^d$ (see \cite{MS19}). We can only prove the $d=2,3$ scaling limit for the domains $\car_L$ because we use the results of \cite{CDH19} and the estimates in \cite{MS19}, which are only known for these domains. The law of $\ol{\phi}(0)$ is not tight as $L \to \infty$ in $d \geq 4$, and so \tref{scal_limit_pt} fails in these dimensions as can be seen by setting $k=1$ and $x_1 = 0$.
			
			\begin{remark} \label{r.regularity.assumption}
				In the case that the hypercube $(-1,1)^d$ is replaced by bounded smooth domain $\Omega$ with $C^{5,\gamma}$ boundary for some $\gamma \in (0,1)$, the solution $u$ of \eref{regularity.assumption} is in $C^{5,\gamma}(\ol{\Omega})$ whenever $f \in C^{1,\gamma}(\ol{\Omega})$ (\cite{GGS10}, Theorem 2.19). Thus \tref{scal_limit_int} holds for the membrane model on $L\Omega \cap \Z^d$ with this assumption on $f$. For the particular case of $d=2,3$ and square or cubic domains, we can use assumption \eref{A2} as a replacement to $u \in C^5([-1,1]^d)$ due to estimates available from \cite{MS19}; however, \tref{scal_limit_pt} largely supersedes \tref{scal_limit_int} in this context.
			\end{remark}
		
			\begin{remark} \label{r.potential.assumption}
				Our arguments are almost entirely quantitative, and implicit constants $C$ in convergence rates are polynomials in the relevant semi-norms of $V$. One could therefore attain results for more general $V$, with diminished convergence rates, by conducting an $L$-dependent cut-off or smoothing to the potential $V$, applying our arguments to the membrane model with regularized $V$, and proving that statistics of interest do not change significantly upon modifying $V$.
			\end{remark}
		\end{subsection}

		\begin{subsection}{Notation and Proof Outline}
			Throughout the article, we work on a cube $\car_L= [-L,L]^d \cap \Z^d$ assumed to be large. We write $x \sim y$ if the points $x,y \in \Z^d$ are adjacent, i.e.\ $|x-y| = 1$ for the $\ell^1$ norm. The (outer) boundary $\pa \car_L$ of the cube $\car_L$ is the set of all $z \in \Z^d \setminus \car_L$ adjacent to an element of $\car_L$ in the integer lattice. We define $\cl_1 \car_L := \car_L \cup \pa \car_L$, and we let $\rho_x$ denote the $\ell^1$ (graph) distance of $x \in \car_L$ to $\pa \car_L$. The second outer boundary $\pa^2 \car_L$ consists of all points in $\Z^d \setminus \car_L$ within graph distance $2$ of a point in $\car_L$. For $U$ a finite subset of $\Z^d$, we denote by $L^2(U)$ the space $\R^U$ with the standard inner product $\inner{\cdot,\cdot}_U$. Expectations with respect to a probability measure $\mu$ are denoted by $\Esp_{\mu}$. For $x \in \Z^d$, we let $\1_x \in \R^{\Z^d}$ be the field with $\1_x(y) = 1$ if $y = x$ and $\1_x(y) = 0$ if $y \ne x$.

			When it is clear, we will make no distinction between an element of $\R^{\Z^d}$ supported on $\car_L$ and its restriction to $\R^{\car_L}$. We will also sometimes restrict an element of $\R^{\cl_1 \car_L}$ to $\R^{\car_L}$ without writing the restriction explicitly. We often free the notation of parameters that stay fixed throughout a section; most commonly, the $L$ dependence will be omitted and the dependence on a parameter $b \in \R^{\car_L}$ (appearing later) will be omitted.
			
			We will approach the study of $\MM_L$ through the distribution of 
			$$
			\Delta \phi(x) := \sum_{y \sim x} (\phi(y) - \phi(x))
			$$
			where the sum is over all $y$ adjacent to $x$ in $\Z^d$. We note that our normalization of $\Delta$ differs from some papers, including \cite{CDH19}.
			More precisely, we let $\Delta_L$, as a map $\R^{\car_L} \to \R^{\car_L}$, be the restriction of $\Delta$ to functions which are $0$ outside of $\car_L$. The map $\Delta_L$ is a linear bijection: the inverse map $\Delta_L^{-1}$ is given by $\Delta_L^{-1} \eta = \phi$ where $\phi$ is the unique solution of
			\begin{equation} \label{e.phidir} \begin{cases} 
					\Delta \phi(x) = \eta(x), \quad &x \in \car_L,\\
					\phi(x) = 0, \quad &x \in \partial \car_L.
			\end{cases} \end{equation}
			We also define for each $z \in \partial \car_L$ the map $\alpha_L^z(\eta) = \Delta (\Delta_L^{-1} \eta)(z)$. That is, we take $\phi = \Delta_L^{-1} \eta$ for any $\eta \in \R^{\car_L}$, extend $\phi$ by $0$ outside $\car_L$, and let $\alpha^z_L(\eta) = \Delta \phi(z)$. We call $\eta$ the Laplacian field associated to $\phi$. We sometimes refer to $\eta(x)$ as a "spin".
			
			For $z \in \pa \car_L$, we let $P^z_L \in \R^{\cl_1 \car_L}$ be the Poisson kernel, which is discrete harmonic in $\car_L$ and equal to $\1_z$ on $\pa \car_L$. We can compute $\alpha^z_L$ in coordinates as $\alpha^z_L(\eta) = -\inner{P^z_L, \eta}_{\car_L}$. Indeed, we have $P^z_L(x) = \Gamma_L(\tilde{z}, x)$ for $x \in \car_L$, where $-\Gamma_L$ is the Dirichlet Green's function of $\Delta$ on $\car_L$ and $\tilde{z} \in \car_L$ is the unique interior point adjacent to $z$. This is because $x \mapsto \Gamma_L(\tilde{z},x) + \1_z(x)$ is discrete harmonic in $\car_L$ and has the same boundary values as $P^z_L$. It follows that
			$$
			\alpha^z_L(\eta) = -\Delta \(\sum_{x \in \car_L} \Gamma_L(\cdot, x)\eta(x)\)(z) = -\sum_{x \in \car_L} \Gamma_L(\tilde{z},x) \eta(x) = -\inner{P^z_L, \eta}_{\car_L}.
			$$
			
			With these definitions in place, we pushforward $\mathrm{MM}_{L}$ by $\Delta_L$ to get on $\R^{\car_L}$ the measure
			\begin{equation} \label{e.Qdef}
				\Q_{L}(d\eta) := \frac{1}{\K(\car_L)} \exp\(-\sum_{x \in \car_L} V(\eta(x)) - \sum_{z \in \pa \car_L} V(\alpha^z_L(\eta))\) \prod_{x \in \car_L} d\eta(x)
			\end{equation}
			for the partition function $\K(\car_L) =|\det \Delta_L|  \mathrm{Z}(\car_L)$. We also introduce the tilted measure $\Q^b_L(d\eta)$ for any $b \in \R^{\car_L}$, defined by
			\begin{equation} \label{e.Qbdef}
				\Q^b_L(d\eta) := \frac{\K(\car_L)}{\K(\car_L,b)}  \exp\(\sum_{x \in \car_L} b(x) \eta(x)\) \Q_L(d\eta),
			\end{equation}
			where $\K(\car_L,b)$ is a constant making $\Q^b_L$ into a probability measure. The central object of study is the cumulant generating function (c.g.f.) of $\Q_L$, which is
			$$
			a \in \R^{\car_L} \mapsto \log \Esp_{\Q_L} \exp \( \sum_{x \in \car_L} a(x) \eta(x) \) = \log \frac{\K(\car_L, a)}{\K(\car_L)}.
			$$
			Note that for any $a' \in \R^{\car_L}$ we have
			$$
			\sum_{x \in \car_L} a'(x) \phi(x) = \sum_{y \in \car_L} \( \sum_{x \in \car_L} -\Gamma_L(y,x) a'(x) \) \eta(y),
			$$
			where $\phi,\eta$ satisfy \eref{phidir}. If $\phi$ is distributed by $\MM_L$, then $\eta$ will be distributed by $\Q_L$, and vice-versa. It follows that the c.g.f. of the membrane model $\MM_L$ evaluated at some $a' \in \R^{\car_L}$ is the c.g.f. of $\Q_L$ evaluated at $a = \Delta_L^{-1} a'$. We will study $\MM_L$ almost entirely through $\Q_L$.
			
			\begin{remark}
				The $\nabla \phi$ model is commonly studied through the distribution of its gradient. One benefit is that in the critical dimension $d=2$, the infinite volume limit of the gradient exists whereas that of the field does not, and many interesting observables are simple, local functions of the gradient. However, the gradient is far from an i.i.d. field due to the many linear dependencies satisfied by gradients, namely that they must sum to $0$ over all closed loops in $\Z^d$. In the membrane model, the situation is partly reversed: many interesting observables are global functions of the Laplacian field, and any field on $\car_L$ is a valid Laplacian field. The former fact will cause difficulties by requiring us to understand the global structure of the law of $\Delta \phi$, but the latter fact will facilitate this understanding.
			\end{remark}
			
			For a function $v$ of $\eta \in \R^{\car_L}$, we will denote by $\pa_{\eta(x)}v(\eta)$ the partial derivative of $v$ in the coordinate $\eta(x)$, when it exists. We let $H^1(\Q^b_L)$ be the closure of smooth, compactly supported functions on $\R^{\car_L}$ under the norm $(\Esp_{\Q^b_L} |v(\eta)|^2 )^{1/2} + (\sum_{x \in \car_L} \Esp_{\Q^b_L} |\pa_{\eta(x)}v(\eta)|^2 )^{1/2}$. We also often use $H^1(\Q^b_L; \R^{\car_L})$, which is a vector-valued version of $H^1(\Q^b_L)$. For $v \in H^1(\Q^b_L; \R^{\car_L})$, we consider $v = v(x,\eta)$ as a function of both discrete space $x \in \car_L$ and the Laplacian field $\eta \in \R^{\car_L}$.
			
			In \sref{HS_sec}, we will derive the Helffer-Sj\"ostrand representation of the cumulant generating function.
			$$
			\log \frac{\K(\car_L, a)}{\K(\car_L)} = \int_0^1 (1-r) \Var_{\Q^{ra}_L} \left [ \sum_{x \in \car_L} a(x) \eta(x) \right] dr = \int_0^1 (1-r)\inner{a, \Esp_{\Q^{ra}_L} \mathrm{HS}^{ra} a}_{\car_L} dr
			$$
			where $\mathrm{HS}^{b}$ is an operator $\R^{\car_L} \to H^1(\Q^b_L; \R^{\car_L})$. We think of $\mathrm{HS}^b$ as an operator dependent on a random environment $\eta$, and averaging over the randomness $\eta \sim \Q^b_L$ gives $\Esp_{\Q^b_L} \mathrm{HS}^b$, which is the covariance operator $\R^{\car_L} \to \R^{\car_L}$ of $(\eta(x))_{x \in \car_L}$. In the Gaussian case, the random environment and the $b$ dependence is not present for $\mathrm{HS}^b$, which is itself the covariance matrix.

			The operator $\mathrm{HS}^b$ has a variational characterization. For $a,b \in \R^{\car_L}$, define $\mathcal{E}_a(\cdot; \Q^b_L)$ on $H^1(\Q^b_L; \R^{\car_L})$ by
			\begin{equation} \begin{aligned} \label{e.E0overview}
					\mathcal{E}_a(v; \Q^b_L)  &= \frac12 \Esp_{\Q^b_L} \left[\sum_{x,y \in \car_L} | \pa_{\eta(x)} v(y, \eta)|^2\right] \\ &\quad + \frac12 \Esp_{\Q^b_L} \left[ \sum_{x \in \car_L} V''(\eta(x)) |v(x,\eta)|^2 + \sum_{z \in \pa \car_L} V''(\alpha^z_L(\eta)) | \alpha^z_L(v(\cdot,\eta)) |^2 \right] \\
					&\quad -  \Esp_{\Q^b_L} \left[ \sum_{x \in \car_L} a(x) v(x,\eta) \right].
			\end{aligned} \end{equation}
			The energy $\mathcal{E}_a$ is strictly convex and has a unique minimizer, which we define as $\mathrm{HS}^b a$. Because $\mathcal{E}_a$ is quadratic, the map $\mathrm{HS}^b$ is a linear operator.
			
			We will think of each of the four sums within \eref{E0overview} as playing a different role. The first sum, the smoothing term, causes the minimizer $\mathrm{HS}^b a$ of $\mathcal{E}_a(\cdot;\Q^b_L)$ to depend on the global shape of $\Q^b_L$; without this term, we could minimize the terms within the expectation in \eref{E0overview} for each realization of $\eta$ alone. The second sum is the self-interaction term with random positive coefficients $V''(\eta(x))$. The third sum is the boundary term, which introduces interaction. Without this term, the minimizer's values at $(x,\eta)$ would be a function of $\eta(x)$ (and hence we would only need to understand the one-spin marginals of $\Q^b_L$). The last term is the perturbation term, without which the minimizer would be $0$.
			
			The main task is to understand the minimizer $\mathrm{HS}^b a$ of \eref{E0overview}. To do so, we will need some initial understanding of $\Q^b_L$. In \sref{SSlaw_sec}, we will understand the marginal law of a single spin $\eta(x)$ for $x \in \car_L$ far from the boundary. Using this, we then study the minimizer of a simpler energy $\mathcal{F}_a(\cdot;\Q^b_L)$, which is obtained from $\mathcal{E}_a(\cdot;\Q^b_L)$ by deleting the boundary term.
			
			We find in \sref{harm_pert_sec} that if $a$ is harmonic on $\car_L$ (after extending $a$ to $\R^{\cl_1 \car_L}$), then there is an alternative characterization of the minimizer of $\mathcal{E}_a$. We use this to show that $\mathrm{HS}^b a$ is small so long as $\| a \|_{L^2(\pa \car_L)}$ is small.
			
			In \sref{Vpert_sec}, we consider the boundary term in $\mathcal{E}_a(\cdot;\Q^b_L)$ and estimate it in terms of the quantities $\inner{P_L^z,a}_{\car_L}$. By subtracting off an appropriate harmonic function from $a$, we can ensure that $\inner{P_L^z,a}_{\car_L}$ is small enough for us to treat the boundary term as negligible as $L \to \infty$. The harmonic function that we subtract off can be handled by the result in \sref{harm_pert_sec}. Finally, when the boundary term is small, then $\mathrm{HS}^b a$ is approximately equal to the minimizer of the simpler energy without the boundary term, which was understood in \sref{SSlaw_sec}.
			
			In \sref{concl_sec}, we synthesize our results to prove the main theorems. In \sref{append}, we give some basic bounds on the Dirichlet Green's functions for the Laplacian, and we give estimates on the aforementioned decomposition of $a$ into a harmonic piece and a remainder.
			
			We now give some further details about the minimization of $\mathcal{E}_a(\cdot;\Q^b_L)$. For simplicity, consider $b = 0$. First, the quadratic functional $\mathcal{E}_0(\cdot;\Q^0_L)$ is coercive on $H^1(\Q^0_L; \R^{\car_L})$ and is half the second variation of $\mathcal{E}_a(\cdot;\Q^0_L)$. If we can find an ansatz $w$ such that $w$ almost minimizes $\mathcal{E}_a(\cdot;\Q^0_L)$, it necessarily follows that $w$ approximates the true minimizer $\mathrm{HS}^0 a$.
			
			We can construct a natural ansatz $w$ by simply ignoring the boundary term in $\mathcal{E}_a$ and attempting to minimize the rest. For simplicity, we also ignore the smoothing term, though in reality we will have to factor it in. With these conditions, the ansatz is $w(x,\eta) = (V''(\eta(x))^{-1} a(x)$.
			
			What is the typical size of $\inner{P^z_L,w}_{\car_L}$? Assuming no cancellations, the best we can estimate is 
			$$
			|\inner{P^z_L, w}_{\car_L}| \leq \| P^z_L \|_{L^1(\car_L)} \| w \|_{L^\infty(\car_L)} \approx L \| w \|_{L^\infty(\car_L)} \approx L^{-d/2 + 1} \| a \|_{L^2(\car_L)}
			$$
			where in the last approximation we used $\inf V'' > 0$ and assumed $a$ is a "macroscopic field", i.e. $\| a \|_{L^\infty(\car_L)} \approx L^{-d/2} \| a \|_{L^2(\car_L)}$. The boundary term with this guess is
			$$
			\Esp_{\Q^0_L} \sum_{z \in \pa \car_L} V''(\alpha^z_L(\eta)) |\inner{P^z_L, w}_{\car_L}|^2 \leq C \| a \|^2_{L^2(\car_L)} \sum_{z \in \pa \car_L} L^{-d +2} \approx  L\| a \|_{L^2(\car_L)}^2,
			$$
			which is very large, meaning that our ansatz is flawed. It turns out that if $a$ is harmonic and macroscopic, this is essentially the truth in the sense that the boundary term dominates, which explains why we are able to show in \sref{harm_pert_sec} that the minimizer is small in this case.
			
			To avoid this issue, we instead break down $a$ into two pieces: $a = K_L a + K_L^\perp a$. Letting $\mathcal{H}(\car_L)$ denote fields which are discrete harmonic on $\car_L$, the first piece is defined by
			$$
			K_L a = \text{argmin}_{\tilde{a} \in \mathcal{H}(\car_L)} \| \tilde{a} - a \|_{L^2(\cl_1 \car_L)}^2,
			$$
			which is the harmonic Bergman projection, the operator analyzed in \sref{append}. Using that $\mathrm{HS}^0$ is linear, we can minimize $\mathcal{E}_{K_L a}$ and $\mathcal{E}_{K_L^\perp a}$ separately. We handle the $K_L a$ term with \sref{harm_pert_sec}, and so we set $a = K_L^\perp a$ in what follows.
			
			Observe that
			$$
			\inner{K_L^\perp a, P^z_L}_{\car_L} =  \inner{K_L^\perp a, P^z_L}_{\cl_1 \car_L} - K_L^\perp a (z) = - K_L^\perp a (z)
			$$
			for any $z \in \pa \car_L$, which follows from the fact that $K_L^\perp$ is the $L^2(\cl_1(\car_L))$ projection onto the orthogonal complement of $\mathcal{H}({\car_L})$, which is the linear span of $\{P^z_L\}_{ z \in \pa \car_L}$. If $a = K^\perp_L a$ is macroscopic, we thus expect that typically
			$$
			|\inner{a, P^z_L}_{\car_L}| \approx L^{-d/2} \| a \|_{L^2(\car_L)}.
			$$
			Returning to our guess $w$, we have
			$$
			\inner{w, P^z_L}_{\car_L} = \sum_{x \in  \car_L} \( \frac{1}{V''(\eta(x))} \) P^z_L(x) a(x).
			$$
			Since the above sum is over many different sites $x$, we might expect the random coefficients $(V''(\eta(x))^{-1}$ undergo a law of large numbers effect and so can be treated like a constant. We prove such an effect and, together with our understanding of the marginals of $\eta(x)$ in \sref{SSlaw_sec}, we can prove that typically
			$$
			| \inner{w, P^z_L}_{\car_L}| \leq C |\inner{a, P^z_L}_{\car_L}| + \mathrm{Error} \leq C L^{-d/2} \| a \|_{L^2(\car_L)} +\mathrm{Error}
			$$
			as $L \to \infty$, for a negligible term $\mathrm{Error}$, so long as we modify the ansatz $w$ near $\pa \car_L$. Roughly speaking, the boundary term is then
			$$
			\Esp_{\Q^0_L} \sum_{z \in \pa \car_L} V''(\alpha^z_L(\eta)) |\inner{P^z_L, w}_{\car_L}|^2 \leq C \| a \|_{L^2(\car_L)}^2 \sum_{z \in \pa \car_L} L^{-d} \leq C L^{-1}\| a \|_{L^2(\car_L)}^2,
			$$
			which is small. In truth, we get a weaker inequality due to error terms from other sources, but the above computation illustrates broadly the logic of the proof.
		\end{subsection}
	\end{section}
	
	\begin{section}{The Helffer-Sj\"ostrand Equation} \label{s.HS_sec}
		In this section, we define the Helffer-Sj\"ostrand representation for the c.g.f. and prove well-posedness results. We also define and prove basic properties of the energies associated to the representation. To lighten notation, we will drop the $L$ dependence from $\Q^b_L$, $P^z_L$, and $\alpha^z_L$.
		
		Recall the definition of $\Q^b$ from \eref{Qdef} and \eref{Qbdef}. For any $b \in \R^{\car_L}$, we define the operator $\mathcal{L}_{\Q^b}$ by
		$$
		\mathcal{L}_{\Q^b} v = \Delta_{\eta} v + \sum_{x \in \car_L} \pa_{\eta(x)} \log \Q^b \cdot \pa_{\eta(x)} v,
		$$
		where we have also denoted by $\Q^b$ the Lebesgue density of $\Q^b$. The operator
		$$
		\Delta_{\eta} = \sum_{x \in \car_L} \pa_{\eta(x)}^2
		$$
		is the (continuum) Laplacian on $\R^{\car_L}$, not to be mistaken with the discrete Laplacian. The above definition can easily be extended to define $\mathcal{L}_\mu$ for other measures $\mu$ with a positive, smooth Lebesgue density.
		
		We compute
		\begin{equation}
			\mathcal{L}_{\Q^b} v = \Delta_{\eta} v + \sum_{x \in  \car_L} \( -V'(\eta(x)) + \sum_{z \in \pa \car_L} V'(\alpha^z(\eta)) P^z(x) + b(x) \) \pa_{\eta(x)} v.
		\end{equation}
		Note that
		\begin{equation} \label{e.Ladj}
			\mathcal{L}_{\Q^b} = -\sum_{x \in \car_L} \pa^{\ast,b}_{\eta(x)} \pa_{\eta(x)}
		\end{equation}
		where $\pa^{\ast,b}_{\eta(x)}$ is the $L^2(\Q^b)$-adjoint of $\pa_{\eta(x)}$. Another important operator is the commutator of $\pa_{\eta(x)}$ and $\pa_{\eta(y)}^{\ast,b}$, which is given for any $x,y \in \car_L$ by
		\begin{equation} \label{e.commid} \begin{aligned}
				[\pa_{\eta(x)}, \pa_{\eta(y)}^{\ast,b}] &=  \pa_{\eta(x)} \pa_{\eta(y)} \mathcal{H}^b(\eta) \\ &= V''(\eta(x)) \1_{x=y} + \sum_{z \in \pa \car_L} V''(\alpha^z(\eta)) P^z(x) P^z(y).
		\end{aligned} \end{equation}
		Here $\mathcal{H}^b$ is the Hamiltonian associated to $\Q^b \propto e^{-\mathcal{H}_b}$; see \eref{Hbeq} below.
		
		We pause to state the Bakry-Emery criterion, which will be used throughout the paper.
		\begin{prop}[See e.g.\ \cite{BL00}] \label{p.poincarebe}
			Let $\mu$ be proportional to $\exp(-W(\eta))\prod_{x \in \car_L} d\eta(x)$, for some $W \in C^2(\R^{\car_L})$. Suppose we can find $\lambda > 0$ such that $\pa^2_\eta W \geq \lambda\cdot \mathrm{Id}$, where $\pa^2_\eta W$ denotes the Hessian of $W$. Then $\mu$ satisfies the Poincar\'e inequality with constant $\lambda^{-1}$ and the log-Sobolev inequality with constant $2\lambda^{-1}$. That is (for $\pa$ the gradient in $\R^{\car_L}$)
			\begin{equation}
				\Var_{\mu} f \leq \frac{1}{\lambda} \Esp_{\mu} | \pa f |^2
			\end{equation}
			for all $f \in H^1(\mu)$ and
			\begin{equation}
				\Esp_{\mu} (f^2 \log f^2) \leq \frac{2}{\lambda} \Esp_{\mu} | \pa f |^2
			\end{equation}
			for all $f \in H^1(\mu)$ with $\Esp_{\mu} f^2 = 1$.
		\end{prop}
		
		We can apply the criterion to invert $\mathcal{L}_{\Q^b}$.
		\begin{prop} \label{p.gen_sol}
			For any $b \in \R^{\car_L}$, the measure $\Q^b$ satisfies the Poincar\'e inequality with constant $\inf V''$. Furthermore, the equation $\mathcal{L}_{\Q^b} u = f$ admits a unique centered weak solution $u \in H^1(\Q^b)$ for every $f \in L^2(\Q^b)$ with $\Esp_{\Q^b} f = 0$.
		\end{prop}
		\begin{proof}
			The first part is due to the Bakry-Emery criterion. The Hamiltonian $\mathcal{H}^b$ defining $\Q^b$ in \eref{Qdef} is
			\begin{equation} \label{e.Hbeq}
				\mathcal{H}^b(\eta) = \sum_{x \in \car_L} V(\eta(x)) + \sum_{z \in \pa \car_L} V(\alpha^z(\eta)) - \sum_{x \in \car_L} b(x) \eta(x)
			\end{equation}
			which has second derivative
			$$
			\pa_{\eta(x)}\pa_{\eta(y)} \mathcal{H}^b(\eta) = \1_{x = y} V''(\eta(x)) + \sum_{z \in \car_L} V''(\alpha^z(\eta)) P^z(x)P^z(y).
			$$
			Since $\alpha^z \otimes \alpha^z \geq 0$, we see that $\pa^2_\eta \mathcal{H}^b \geq (\inf V'') \text{Id}$, and so \pref{poincarebe} applies. We can then invert $\mathcal{L}_{\Q^b}$ using the Lax-Milgram theorem. The Poincar\'e inequality provides $H^1(\Q^b)$ coerciveness of the quadratic form $-\Esp_{\Q^b} f \mathcal{L}_{\Q^b}f = \Esp_{\Q^b} | \pa f |^2$ on centered random variables $f$.
		\end{proof}
		
		We give a formal computation that motivates studying $\mathcal{L}_{\Q^b}$ (essentially from \cite{NS97}, Section 1). First, one can show that to compute the c.g.f. of $\Q^b$, it is sufficient to understand variances under the tilted measures $\Q^b$. Then, for a centered $F \in H^1(\Q^b)$, we have
		$$
		\Var_{\Q^b} F = \Esp_{\Q^b} \left[ F (-\mathcal{L}_{\Q^b}) (-\mathcal{L}_{\Q^b})^{-1} F \right] = \sum_{x \in \car_L} \Esp_{\Q^b} \left[ \pa_{\eta(x)} F \cdot \pa_{\eta(x)} (-\mathcal{L}_{\Q^b})^{-1} F \right],
		$$
		where in the last equality we used the representation \eref{Ladj}. Let $u(x,\eta) = \pa_{\eta(x)} (-\mathcal{L}_{\Q^b})^{-1} F$. If we apply $\pa_{\eta(x)}$ to the equation
		$$
		(-\mathcal{L}_{\Q^b}) (-\mathcal{L}_{\Q^b})^{-1} F = F
		$$ 
		and commute the gradient and $\mathcal{L}_{\Q^b}$, we get
		$$
		-\mathcal{L}_{\Q^b} u(x,\eta) - [\pa_{\eta(x)},\mathcal{L}_{\Q^b}] (-\mathcal{L}_{\Q^b})^{-1}F(\eta) = \pa_{\eta(x)} F(\eta).
		$$
		For the commutator above, we compute
		$$
		[\pa_{\eta(x)},-\mathcal{L}_{\Q^b}] = \sum_{y \in \car_L} [\pa_{\eta(x)}, \pa_{\eta(y)}^{\ast,b} \pa_{\eta(y)}] = \sum_{y \in \car_L} [ \pa_{\eta(x)}, \pa_{\eta(y)}^{\ast,b}] \pa_{\eta(y)},
		$$
		and we use \eref{commid} on the second term above to see
		$$
		-\mathcal{L}_{\Q^b} u(x,\eta)  + (\pa^2_{\eta} \mathcal{H}^b \cdot  u)(x,\eta) = \pa_{\eta(x)} F(\eta).
		$$
		This is the Helffer-Sj\"ostrand equation for $u$. It follows that
		$$
		\Var_{\Q^b} F = \Esp_{\Q^b} \left[ \inner{f(\cdot,\eta), u(\cdot,\eta)}_{\car_L} \right]
		$$
		for $f(x,\eta) = \pa_{\eta(x)} F(\eta)$, which is the Helffer-Sj\"ostrand representation.
		
		The next two lemmas make rigorous the above computation in the special case $f(x,\eta) = a(x)$.
		\begin{lem} \label{l.HSsol_lem}
			For any $a,b \in \R^{\car_L}$, the Helffer-Sj\"ostrand equation
			\begin{equation} \label{e.HSeq}
				-\mathcal{L}_{\Q^b} v + \pa_\eta^2 \mathcal{H}^b \cdot v = a
			\end{equation}
			has a unique weak solution $\mathrm{HS}^b a = u^b_a \in H^1(\Q^b; \R^{\car_L})$. Also, we have $u^b_a(x,\eta) = \pa_{\eta(x)}v^b_a(\eta)$, where $v^b_a$ solves (weakly)
			\begin{equation} \label{e.vba_eq}
				-\mathcal{L}_{\Q^b} v^b_a  = \mathrm{Const} + \sum_{x \in \car_L} a(x) \eta(x).
			\end{equation}
			The constant $\mathrm{Const}$ is chosen so that the RHS above has expectation $0$ under $\Q^b$.
		\end{lem}
		\begin{proof}
			Uniqueness and existence of $u^b_a$ follows from $\pa_\eta^2 \mathcal{H}^b\geq (\inf V'') \text{Id}_{\car_L \times \car_L}$ a.s.\ in $\eta$ and from Lax-Milgram. Given $v^b_a \in H^1(\Q^b)$ solving \eref{vba_eq} (such a solution exists by \pref{gen_sol}), we will show that $\pa_{\eta(x)} v^b_a$ satisfies equation \eref{HSeq} in a weak sense, i.e.\ integrated against $H^1(\Q^b; \R^{\car_L})$ functions. From the weak form of the equation for $v^b_a$, we have
			$$
			\sum_{x \in \car_L} \Esp_{\Q^b} \left[ \pa_{\eta(x)} v^b_a \cdot \pa_{\eta(x)} \varphi \right] = \mathrm{Const} \cdot \Esp_{\Q^b} \varphi + \sum_{x \in \car_L} a(x) \Esp_{\Q^b} \left[ \eta(x) \varphi \right]
			$$
			for all $\varphi \in H^1(\Q^b)$. Choose $\varphi = \pa^{\ast,b}_{\eta(x_0)} \psi$ for $\psi \in C^\infty_c(\R^{\car_L})$ to see
			\begin{align*}
				\lefteqn{\Esp_{\Q^b} \left[ \pa_{\eta(x)} v^b_a \cdot \pa_{\eta(x)} \varphi \right]} \quad & \\ &=  \Esp_{\Q^b} \left[ \pa_{\eta(x)} v^b_a \cdot \(\pa^{\ast,b}_{\eta(x_0)} \pa_{\eta(x)} + [\pa_{\eta(x)}, \pa^{\ast,b}_{\eta(x_0)}]\) \psi \right] \\
				&= \Esp_{\Q^b} \left[ \pa_{\eta(x_0)} v^b_a \cdot \pa^{\ast,b}_{\eta(x)} \pa_{\eta(x)} \psi \right] + \Esp_{\Q^b}\left[ \pa_{\eta(x)} v^b_a \cdot [\pa_{\eta(x)}, \pa^{\ast,b}_{\eta(x_0)}] \psi \right].
			\end{align*}
			In the last line we moved the $\pa_{\eta(x)}$ derivative onto $\psi$, used $[\pa^{\ast,b}_{\eta(x)}, \pa^{\ast,b}_{\eta(x_0)}] = 0$, and then put the derivative $\pa_{\eta(x_0)}$ onto $v^b_a$. We then use \eref{commid} and sum over $x$ to see
			$$
			\Esp_{\Q^b} \left[ \pa_{\eta(x_0)} v^b_a (-\mathcal{L}_{\Q^b} \psi) \right] + \Esp \left[ \sum_{x \in \car_L} \(\pa^2_{\eta(x),\eta(x_0)} \mathcal{H}^b\) (\pa_{\eta(x)} v^b_a) \cdot \psi \right] = a(x_0) \Esp_{\Q^b} \left[ \psi \right].
			$$
			Let $\tilde{u}(x_0, \eta) := \pa_{\eta(x_0)} v^b_a(\eta)$. We now let $\psi = \psi_{x_0}$ depend on $x_0$ and sum over $x_0$. Note that the Hessian of $\mathcal{H}^b$ is self-adjoint, and so we get
			\begin{equation} \label{e.tildeu.weak}
			\Esp_{\Q^b} \inner{\tilde{u}(\cdot,\eta), (-\mathcal{L}_{\Q^b} + \pa^2_{\eta} \mathcal{H}^b)  \psi_{\cdot}(\eta) }_{\car_L} = \Esp_{\Q^b} \inner{a, \psi_{\cdot}}_{\car_L}.
			\end{equation}
			This is a weak form of \eref{HSeq} which also implies that $\tilde u$ is in $H^1(\Q^b;\R^{\car_L})$. Indeed, the LHS of \eref{tildeu.weak} as a function of $\psi$ defines a continuous linear functional on a dense subspace of $H^1(\Q^b; \R^{\car_L})$, and so it can be uniquely extended to the full space. Since $\pa^2_\eta \mathcal{H}^b \geq (\inf V'') \mathrm{Id}$, we have
			$$
			\| \tilde u \|_{H^1(\Q^b;\R^{\car_L})} \leq C \sup_{ \| \psi \|_{H^1(\Q^b; \R^{\car_L})} \leq 1}\Esp_{\Q^b} \inner{\tilde{u}(\cdot,\eta), (-\mathcal{L}_{\Q^b} + \pa^2_{\eta} \mathcal{H}^b)  \psi_{\cdot}(\eta) }_{\car_L},
			$$
			which is finite by \eref{tildeu.weak}. We can then integrate by parts in \eref{tildeu.weak} to show that $\tilde{u}$ solves \eref{HSeq} weakly in $H^1(\Q^b; \R^{\car_L})$, and by uniqueness we have $u^b_a = \tilde{u}$.
		\end{proof}
		
		We will write the solution of \eref{HSeq} as $u^b_a$ or $\mathrm{HS}^b a$. When $b$ is understood, or we are considering solutions of \eref{HSeq} where the operator $\mathcal{L}_{\Q^b}$ is replaced by $\mathcal{L}_\mu$ for some reference measure $\mu$, we will often just write $u_a$.
		
		\begin{lem}
			We can express the cumulant generating function of $\Q^0$ as
			\begin{equation} \label{e.cumulant_var}
				\log \frac{\K(\car_L, a)}{\K(\car_L, 0)} = \int_0^1 (1 - r) \Var_{\Q^{ra}} \left[\sum_{x \in \car_L} a(x) \eta(x) \right] dr.
			\end{equation}
			The variances can be expressed via the Helffer-Sj\"ostrand representation by
			\begin{equation} \label{e.var_form}
				\Var_{\Q^{ra}} \left[\sum_{x \in \car_L} a(x) \eta(x) \right] = \inner{a, \Esp_{\Q^{ra}} u^{ra}_a}_{\car_L} =  \inner{a, \Esp_{\Q^{ra}} \mathrm{HS}^{ra} a}_{\car_L} .
			\end{equation}
		\end{lem}
		\begin{proof}
			We start by computing
			$$
			\frac{d}{dr} \log \K(\car_L, ra) = \Esp_{\Q^{ra}} \left[\sum_{x \in \car_L} a(x) \eta(x) \right].
			$$
			Since $\Q^0$ is invariant under $\eta \mapsto -\eta$, we see $\frac{d}{dr} \log \K(\car_L, ra) = 0$ at $r = 0$. We can take another derivative in $r$ to see
			$$
			\(\frac{d}{dr}\)^2  \log \K(\car_L, ra) = \Var_{\Q^{ra}} \left[\sum_{x \in \car_L} a(x) \eta(x) \right].
			$$
			We conclude \eref{cumulant_var} by integrating in $r$.
			
			We now prove the variance formula \eref{var_form}. Let $F(\eta) = \sum_{x \in \car_L} a(x) \eta(x)$. We use the formula
			$$
			u^{ra}_a(x,\cdot) = \pa_{\eta(x)} (-\mathcal{L}_{\Q^{ra}})^{-1} \( F - \Esp_{\Q^{ra}} F \)
			$$
			which was proved in \lref{HSsol_lem}.  Since $\mathcal{L}_{\Q^{ra}}$ is $L^2(\Q^{ra})$ self-adjoint and the representation \eref{Ladj} holds, we have
			\begin{equation*} \begin{aligned}
					\Var_{\Q^{ra}} F &= \Esp_{\Q^{ra}} \left[ F(\eta)(F(\eta) - \Esp_{\Q^{ra}} F)\right] \\ &= -\Esp_{\Q^{ra}} \left[ \mathcal{L}_{\Q^{ra}} F(\eta)(-\mathcal{L}_{\Q^{ra}})^{-1}(F(\eta) - \Esp_{\Q^{ra}} F)\right] \\
					&= \Esp_{\Q^{ra}} \left[ \sum_{x \in \car_L}\( \pa_{\eta(x)} F(\eta) \) u^{ra}_a(x,\eta) \right] = \Esp_{\Q^{ra}} \inner{a, u^{ra}_a(x,\eta)}_{\car_L}.
		\end{aligned} \end{equation*} \end{proof}
		
		We will now define three energies and state some basic theorems about their minimizers. Two of the energies will depend on a background measure $\Q^b$ for some $b \in \R^{\car_L}$, and we will sometimes omit this dependence. As in \eref{E0overview}, for any $a \in \R^{\car_L}$, the energy $\mathcal{E}_a(\cdot;\Q^b) : H^1(\Q^b; \R^{\car_L}) \to \R$ is defined by
		\begin{equation} \begin{aligned}\label{e.Eadef}
				\mathcal{E}_a(v; \Q^b)  &= \frac12 \Esp_{\Q^b} \left[ \sum_{x, y \in \car_L} |\pa_{\eta(y)} v(x,\eta)|^2 +  \sum_{x \in \car_L} V''(\eta(x)) |v(x,\eta)|^2 \right] \\ &\quad+ \frac12 \Esp_{\Q^b} \left[ \sum_{z \in \pa \car_L} V''(\alpha^z(\eta)) |\alpha^z_L(v(\cdot,\eta))|^2 \right] - \Esp_{\Q^b} \left[ \sum_{x \in \car_L} a(x) v(x,\eta) \right].
		\end{aligned} \end{equation}
		The minimizer of $\mathcal{E}_a$ will be denoted by $u_a$ (we will see that it is the solution to the Helffer-Sj\"ostrand equation momentarily).
		
		We also introduce a version of $\mathcal{E}_a$ without the boundary term:
		\begin{equation} \label{e.F0def}
			\mathcal{F}_{0}(v; \Q^b) = \frac12 \Esp_{\Q^b} \left[ \sum_{x, y \in \car_L} |\pa_{\eta(y)} v(x,\eta)|^2 +  \sum_{x \in \car_L} V''(\eta(x)) |v(x,\eta)|^2 \right]
		\end{equation}
		and
		\begin{equation} \label{e.Fadef}
			\mathcal{F}_{a}(v; \Q^b) = \mathcal{F}_{0}(v; \Q^b) - \Esp_{\Q^b} \left[ \sum_{x \in \car_L} a(x) v(x,\eta) \right].
		\end{equation}
		We denote by $\ol{u}_{a}$ the minimizer of $\mathcal{F}_a(\cdot;\Q^b)$.
		
		\begin{prop} \label{p.energy_bds}
			The functionals $\mathcal{E}_{a}(\cdot;\Q^b)$ and $\mathcal{F}_{a}(\cdot;\Q^b)$ are strictly convex and continuous on $H^1(\Q^b; \R^{\car_L})$. We have the lower bounds
			\begin{equation} \label{e.lowbd1}
				\inf_{v \in H^1(\Q^b; \R^{\car_L})} \mathcal{E}_{a}(v; \Q^b) \geq \inf_{v \in H^1(\Q^b; \R^{\car_L})} \mathcal{F}_{a}(v; \Q^b) \geq -\frac{1}{2\inf V''}\| a \|^2_{L^2(\car_L)}.
			\end{equation}
			The minimizers $u_{a}$ and $\ol{u}_{a}$ are weak solutions of the following equations:
			\begin{align} \label{e.HSeqs.minimizer1}
				&-\mathcal{L}_{\Q^b} u_{a}(x,\eta) + (\pa_\eta^2 \mathcal{H}^b \cdot u_a) (x,\eta) = a(x) \\ 
				&-\mathcal{L}_{\Q^b} \ol{u}_{a}(x,\eta) + V''(\eta(x)) \ol{u}_{a}(x,\eta) = a(x).
			\end{align}
			For any $v \in H^1(\Q^b; \R^{\car_L})$, we have the quadratic response identities
			\begin{equation}
				\mathcal{E}_{a}(v ;\Q^b) - \mathcal{E}_{a}(u_{a};\Q^b) = \mathcal{E}_{0}(v - u_{a};\Q^b), \quad \mathcal{F}_{a}(v; \Q^b) - \mathcal{F}_{a}(\ol{u}_{a};\Q^b) = \mathcal{F}_{0}(v - \ol{u}_{a};\Q^b).
			\end{equation}
			Consequently (letting $v = 0$), we have
			\begin{equation} \label{e.upbd1}
				\mathcal{E}_{0}(u_a;\Q^b) \leq \mathcal{F}_0(\ol{u}_{a};\Q^b) \leq \frac{1}{2\inf V''}\| a \|^2_{L^2(\car_L)}.
			\end{equation}
			Furthermore, $\ol{u}_{a}(x,\eta)$ is between $(\inf V'')^{-1}a(x)$ and $(\sup V'')^{-1}a(x)$ almost surely (in $\eta$) for all $x \in \car_L$.
		\end{prop}
		\begin{proof}
			Continuity and strict convexity of the functionals follow easily from $0 < \inf V'' \leq \sup V'' < \infty$.
			
			The inequalities in \eref{lowbd1} follow from ignoring all terms in the energy except for
			$$
			\Esp_{\Q^b} \left[ \frac12\sum_{x \in \car_L} V''(\eta(x)) |v(x,\eta)|^2 - \sum_{x \in \car_L} a(x) v(x,\eta) \right].
			$$
			The equations for $u_a$ and $\ol{u}_a$ follow from standard first variation arguments. The quadratic response identities follow immediately from expanding the square within the energies and using the equations for $u_a$ and $\ol{u}_a$. The bound \eref{upbd1} follows from plugging $v = 0$ into the quadratic response identities and using \eref{lowbd1}.
			
			Finally, we prove the bound on $\ol{u}_{a}$. Let $[c_x,d_x]$ be the interval between $(\inf V'')^{-1}a(x)$ and $(\sup V'')^{-1}a(x)$. Define $\chi_x : \R \to \R$ by $\chi_x(y) = \max(c_x, \min(y,d_x))$. Then $(x,\eta) \mapsto \chi_x(\ol{u}_{a}(x,\eta))$ is a minimizer of $\mathcal{F}_{a}$. This is because $\pa_{\eta(x)} \chi_x(\ol{u}_{a}(x,\eta)) = 0$ almost surely when $\chi(\ol{u}_a)$ is $c$ or $d$, and
			\begin{equation*}\begin{aligned}
					\frac12 \sum_{x \in \car_L} V''(\eta(x))|\chi_x(\ol{u}_a(x,\eta))|^2 - \sum_{x \in \car_L} a(x) \chi_x(\ol{u}_a(x,\eta)) &\leq \frac12 \sum_{x \in \car_L} V''(\eta(x))|\ol{u}_a(x,\eta)|^2 \\ &\quad - \sum_{x \in \car_L} a(x) \ol{u}_a(x,\eta).
			\end{aligned}\end{equation*}
			By uniqueness of the minimizer, we then have $\chi_x(\ol{u}_a(x,\eta)) = \ol{u}_{a}(x,\eta)$ a.s.
		\end{proof}
		
		The main importance of $\ol{u}_a$ is that it is significantly easier to understand than $u_a$, and so it serves as a building block for constructing near-minimizers of $\mathcal{E}_a$. In particular, we will approximate $\ol{u}_a$ in terms of the following 1D object. For any $\beta \in \R$, we let $\nu^\beta$ be the probability measure on $\R$ proportional to $\exp(-V(\xi) + \beta \xi) d\xi$. Then, for $v \in H^1(\nu^{\beta})$ we let
		\begin{equation} \label{e.G0def}
			\mathcal{G}_0(v; \nu^\beta) = \Esp_{\nu^\beta} \left[\frac12 |\pa_\xi v|^2 + \frac12 V''(\xi) |v(\xi)|^2 \right]
		\end{equation}
		and, for any $\gamma \in \R$, we let $\mathcal{G}_\gamma(v;\nu^\beta) = \mathcal{G}_0(v; \nu^\beta) - \gamma \Esp_{\nu^\beta} v$. We let $U^\beta$ be the minimizer of $\mathcal{G}_1(\cdot; \nu^\beta)$.
		
		\begin{prop} \label{p.Genergy}
			The energy $\mathcal{G}_1(\cdot; \nu^\beta)$ is continuous on $H^1(\nu^\beta)$ and strictly convex. The associated minimizer has $(\sup V'')^{-1} \leq U^\beta \leq (\inf V'')^{-1}$ a.s. and $\| \pa_\xi U^\beta \|_{L^\infty(\R)} \leq C$.
		\end{prop}
		\begin{proof}
			The proof of $(\sup V'')^{-1} \leq U^\beta \leq (\inf V'')^{-1}$ is exactly the same as the corresponding bound on $\ol{u}_a$ in \pref{energy_bds}. To prove the bound on the gradient of $U^\beta$, consider the dynamics
			$$
			d\xi_t = (- V'(\xi_t) + \beta)dt + \sqrt{2}dB_t
			$$
			for a standard Brownian motion $B_t$. The operator $-\mathcal{L}_{\nu^\beta} = -\Delta_\xi + (V'(\xi) - \beta) \pa_\xi$ is the generator of the dynamics. As $U^\beta$ solves $(-\mathcal{L}_{\nu^\beta} + V'')U^\beta = 1$, it enjoys the Feynman-Kac representation
			$$
			U^\beta(\xi) = \Esp^{\xi_0 = \xi} \left[ \int_0^\infty \exp\(-\int_0^t V''(\xi_s) ds\) dt \right].
			$$ 
			And so
			$$
			\pa_{\xi} U^\beta(\xi) = -\Esp^{\xi_0 = \xi} \left[ \int_0^\infty \( \exp\(-\int_0^t V''(\xi_s) ds \) \int_0^t V^{(3)}(\xi_s) J_{0,s}(\xi) ds \) dt \right]
			$$
			where $J_{0,s}(\xi)$ is the "derivative" flow of $\xi_t$, which solves
			\begin{equation*} \begin{cases}
					\pa_t J_{0,t}(\xi) = -V''(\xi_t)J_{0,t}(\xi), \\
					J_{0,0}(\xi) = 1.
			\end{cases} \end{equation*}
			It is clear that $J_{0,t}(\xi)$ is bounded uniformly in $(\xi_t)_{t \geq 0}$ and $t$ since its equation is contractive, and $V^{(3)}$ is bounded by assumption. Thus
			$$
			|\pa_{\xi} U^\beta(\xi)| \leq C \Esp^{\xi_0 = \xi}\left[  \int_0^\infty  \exp\(-\int_0^t V''(\xi_s) ds  \) tdt \right] \leq C\int_0^\infty t e^{- \inf V'' t} dt \leq C.
			$$
		\end{proof}
	\end{section}
	
	\begin{section}{Minimization for Smooth Harmonic Perturbations} \label{s.harm_pert_sec}
		In this section, we let $a = \sum_{z \in \pa \car_L} \gamma_z P^z$ for some coefficients $(\gamma_z)_{z \in \pa \car_L}$ and the Poisson kernel $P^z = P^z_L$ on $\car_L$, and we want to prove that the minimizer $u_a$ of $\mathcal{E}_a(\cdot;\Q^b)$ is small. This will prove $\mathrm{HS}^b a$ is small and ultimately allow us to discard harmonic components of more general $a$ in our calculations. The result will be used in \sref{concl_sec} to prove the main theorems, but is independent from \sref{SSlaw_sec} or \sref{Vpert_sec}. We omit the dependence on $\Q^b$ below, and all constants $C$ are independent of $b$.
		
		We first reduce the problem to the same problem for a Gaussian membrane model. To that end, let $\mathcal{E}'_a$ be the energy
		\begin{equation} \begin{aligned}
				\mathcal{E}'_a(v) &= \frac12 \Esp_{\Q^b} \left [  \sum_{x \in \car_L} (\inf V'') |v(x,\eta)|^2 + \sum_{z \in \pa \car_L} (\inf V'') | \alpha^z(v(\cdot,\eta))|^2 \right ] \\ &\quad- \Esp_{\Q^b} \left [ \sum_{x \in \car_L} a(x) v(x,\eta) \right ].
		\end{aligned} \end{equation}
		Note that $\mathcal{E}'_0(v)$, and its perturbed versions $\mathcal{E}'_a$, correspond to the energies $\mathcal{E}_a$ with the quadratic potential $\xi \mapsto (\inf V'')|\xi|^2$ in place of $V$. We do not include the smoothing term $|\pa_{\eta} v|^2$ because all minimizers are independent of $\eta$ (it would also be of no consequence to include it). We have
		\begin{equation} \label{e.EE'}
			\mathcal{E}_a(v) \geq  \mathcal{E}'_a(v).
		\end{equation}
		Let $u_a$ be the minimizer of $\mathcal{E}_a$ and let $u'_a$ be the minimizer of $\mathcal{E}'_a$.
		
		There is a small technical issue in that the representation $a = \sum_{z \in \pa \car_L} \gamma_z P^z$ is not unique given just the data $a \in \R^{\car_L}$. It can happen that $P^z = P^w$ in $\car_L$ for $z \ne w$, for example if $z$ and $w$ are adjacent to the same corner of the square $\car_{L}$. The below proposition will be valid for any choice of extension of $a$ to $\cl_1 \car_L$.
		\begin{prop} \label{p.Estiffprop}
			Let $a = \sum_{z \in \pa \car_L} \gamma_z P^z$ and $u_a = \mathrm{HS}^b a$. Then there is a constant $C$ such that
			\begin{equation} \label{e.Estiff1}
				\mathcal{E}_0(u_a; \Q^b) \leq C \| a \|_{L^2(\car_L)} \| \gamma \|_{L^2(\pa \car_L)}.
			\end{equation}
		\end{prop}
		
		Note that the estimate \eref{Estiff1} has the boundary term $\| \gamma \|_{L^2(\pa \car_L)}$, which should be much smaller than $\| a \|_{L^2(\car_L)}$ provided $a$ is macroscopic or smooth in a sense. This makes the estimate much stronger than the trivial estimate \eref{upbd1} with RHS $C \| a \|_{L^2(\car_L)}^2$.
		
		\begin{proof}
			It is enough to provide a lower bound for $\mathcal{E}_a$, since we have $\mathcal{E}_0(u_a) = -\mathcal{E}_a(u_a)$ by comparison with $v=0$ and quadratic response. By \eref{EE'}, it is enough to bound $\mathcal{E}'_a(u'_a)$. The proof proceeds by rewriting the minimization problem in terms of the boundary values $\gamma_z$ of $a$. 
			
			Define a strictly convex, continuous functional $\mathcal{B}_\gamma : \R^{\pa \car_L} \to \R$ by
			$$
			\mathcal{B}_\gamma(\epsilon) = \frac12 (\inf V'') \bigg \| \sum_{z \in \pa \car_L} \epsilon_z P^z  \bigg \|_{L^2(\cl_1 \car_L)}^2 - \sum_{z \in \pa \car_L} \epsilon_z \gamma_z.
			$$
			Let $\lambda \in \R^{\pa \car_L}$ be the minimizer of $\mathcal{B}_\gamma$. We claim that $\ol{\lambda} = \sum_{z \in \pa \car_L} \lambda_z P^z \in \R^{\car_L}$ is equal to the minimizer $u'_a$ of $\mathcal{E}'_a$ within $\car_L$. Indeed, the critical point $\lambda$ solves
			$$
			\sum_{w \in \pa \car_L} \lambda_w \inner{P^z, P^w}_{\cl_1 \car_L} = \frac{\gamma_z}{\inf V''} \quad \forall z \in \pa \car_L,
			$$
			which we rewrite to
			\begin{equation} \label{e.Epvar}
				\lambda_z + \sum_{w \in \pa \car_L} \lambda_w \inner{P^z, P^w}_{\car_L} = \frac{\gamma_z}{\inf V''} \quad \forall z \in \pa \car_L.
			\end{equation}
			After multiplying \eref{Epvar} by $P^z(x)$ and summing over $z$, we get that
			$$
			\ol{\lambda}(x) + \sum_{z \in \pa \car_L} \inner{P^z, \ol{\lambda}}_{\car_L} P^z(x) = \frac{a(x)}{\inf V''} \quad \forall x \in \car_L.
			$$
			This equation also characterizes minimizers of $\mathcal{E}'_a$, so we conclude $\ol{\lambda} = u'_a$.
			
			Note that (using $\sum_{z \in \pa \car_L} \epsilon_z P^z(x) \geq \sum_{z \in \pa \car_L} \epsilon_z \1_{z}(x)$ for all $x \in \cl_1 \car_L$ and Cauchy-Schwarz)
			$$
			\mathcal{B}_\gamma(\epsilon) \geq \frac12 (\inf V'') \| \epsilon \|_{L^2(\pa \car_L)}^2 -  \| \epsilon \|_{L^2(\pa \car_L)}  \| \gamma \|_{L^2(\pa \car_L)} \geq -\frac{1}{2 (\inf V'')}  \| \gamma \|^2_{L^2(\pa \car_L)}.
			$$
			By quadratic response we have
			$$
			\mathcal{B}_0(\lambda) \leq \mathcal{B}_\gamma(0) - \mathcal{B}_\gamma(\lambda) \leq \frac{1}{2 (\inf V'')}  \| \gamma \|^2_{L^2(\pa \car_L)},
			$$
			and thus
			\begin{equation} \label{e.Hgap}
				\| u'_a \|_{L^2(\car_L)}^2 \leq \frac{2}{\inf V''} \mathcal{B}_0(\lambda) \leq \frac{1}{(\inf V'')^2}\| \gamma \|^2_{L^2(\pa \car_L)}.
			\end{equation}
			Therefore, we have
			$$
			\mathcal{E}'_a(u'_a) \geq -\inner{u'_a, a}_{\car_L} \geq -C \| a \|_{L^2(\car_L)} \| \gamma \|_{L^2(\pa \car_L)},
			$$
			and the claim \eref{Estiff1} is established.
		\end{proof}
	\end{section}
	
	\begin{section}{Marginals and the Simplified Energy} \label{s.SSlaw_sec}
		The goal of this section is to understand the energy $\mathcal{F}_a(\cdot;\Q^b)$ defined in \eref{Fadef}. To do this, we will need to understand the marginal law $\mu^b_{x_0}$ of $\eta(x_0)$ under $\Q^b$ for $x_0 \in \car_L$. Specifically, we will prove an approximation
		\begin{equation} \label{e.nub_def}
			\mu^b_{x_0}(d\eta(x_0)) \approx \nu^{b(x_0)}(d\eta(x_0)) := \frac{1}{\mathrm{J}(b(x_0))} \exp(b(x_0) \eta(x_0) - V(\eta(x_0))) d\eta(x_0)
		\end{equation}
		where $\mathrm{J}(b(x_0))$ is a normalizing constant and $x_0$ is not too close to $\pa \car_L$. This approximation would be exact if the boundary term were not present in the energy $\mathcal{H}^b$ defining $\Q^b$.
		
		We will also show that the "off-diagonal" derivatives $\pa_{\eta(y)} v(x,\eta)$ for $x \ne y$ within the smoothing term of $\mathcal{F}_a(v;\Q^b)$ are not consequential for the minimization problem. After these off-diagonal terms are eliminated and $\mu^b_{x_0} \approx \nu^{b(x_0)}$ is applied, the energy $\mathcal{F}_a(\cdot;\Q^b)$ just becomes a weighted sum of simple energies like $\mathcal{G}_1(\cdot;\nu^{b(x)})$ which can be minimized independently.
		
		Our main tool for both of these results is a special profile $e_{x_0} \in \R^{\car_L}$, which we now describe.
		
		\begin{subsection}{The Special Profile}
			To motivate the construction of the special profile, consider the following computation of the derivative of the marginal of $\eta(x_0)$. We have
			$$
			\mu^b_{x_0}(\eta(x_0)) = \frac{1}{\K(\car_L,b)} \int_{\R^{\car_L \setminus \{x_0\}}} e^{-\mathcal{H}^b(\eta)} \prod_{x \ne x_0} d\eta(x).
			$$
			where we denote by $\mu^b_{x_0}(\eta(x_0))$ the Lebesgue density of the marginal, by abuse of notation. Then for every $\epsilon \in \R$, we have
			$$
			\mu^b_{x_0}(\eta(x_0) + \epsilon) - \mu^b_{x_0}(\eta(x_0)) = \Esp_{\Q^b_c} \left[ e^{\mathcal{H}^b(\eta) - \mathcal{H}^b(\eta + \epsilon \1_{x_0})} - 1 \right] \mu^b_{x_0}(\eta(x_0))
			$$
			where $\Q^b_c$ means $\Q^b$ conditioned on $\eta(x_0)$. And so
			\begin{equation*} \begin{aligned}
					\pa_{\eta(x_0)} \log \mu^b_{x_0}(\eta(x_0)) &= -\Esp_{\Q^b_c} \left[ \pa_{\eta(x_0)} \mathcal{H}^b(\eta) \right] \\ &= - V'(\eta(x_0)) + b(x_0) + \Esp_{\Q^b_c} \left[ \sum_{z \in \pa \car_L} V'(\alpha^z(\eta)) P^z(x_0) \right].
			\end{aligned} \end{equation*}
			The first two terms on the RHS are consistent with our approximation $ \mu^b_{x_0}(\eta(x_0)) \approx \nu^{b(x_0)}$, but the last boundary term needs to be small. We do not know a way to bound this quantity directly to satisfactory accuracy, particularly in $d=2,3$ where the summability properties of $P^z$ are less desirable.
			
			If instead we apply a change of variables $\eta \mapsto \eta - \epsilon e_{x_0}$ for some $e_{x_0} \in \R^{\car_L}$ with $e_{x_0}(x_0) = 0$ to the integral defining $\mu^b_{x_0}(\eta(x_0) + \epsilon)$, we can "move" the derivative off of the boundary term. We get
			$$
			\mu^b_{x_0}(\eta(x_0) + \epsilon) - \mu^b_{x_0}(\eta(x_0)) = \Esp_{\Q^b_c} \left[ e^{\mathcal{H}^b(\eta) - \mathcal{H}^b(\eta + \epsilon \1_{x_0} + \epsilon e_{x_0})} - 1 \right]\mu^b_{x_0}(\eta(x_0)),
			$$
			and so
			\begin{equation*} \begin{aligned}
					\pa_{\eta(x_0)} \log \mu^b_{x_0}(\eta(x_0)) &= -V'(\eta(x_0)) + b(x_0) +  \Esp_{\Q^b_c} \left[ \sum_{x \ne x_0} (-V'(\eta(x)) + b(x)) e_{x_0}(x) \right]\\ &\quad +  \Esp_{\Q^b_c} \left[ \sum_{z \in \pa \car_L} V'(\alpha^z(\eta)) \inner{P^z, \1_{x_0} + e_{x_0}}_{\car_L} \right].
			\end{aligned} \end{equation*}
			If we construct $e_{x_0}$ such that $\inner{\1_{x_0} + e_{x_0}, P^z} = 0$, the boundary term now vanishes. The extra term will be much easier to bound, particularly since $e_{x_0}$ will be small. The above computation suggests that an a priori strong interaction between $\eta(x_0)$ and the $O(L^{d-1})$ boundary terms can be mitigated through a very small variation by $O(L^d)$ other spins $\eta(x)$, $x \ne x_0$. 
			
			\begin{prop} \label{p.berg_con}
				Recall that $\rho_{x_0}$ denotes the distance between $x_0$ and $\pa \car_L$. For each ${x_0 \in \car_L}$ not adjacent to $\pa \car_L$, there exists $e_{x_0} \in \R^{\cl_1 \car_L}$ such that $e_{x_0}(x_0) = 0$, $\| e_{x_0} \|_{L^2(\cl_1 \car_L)} \leq C\rho_{x_0}^{-d/2}$ and 
				$$
				\inner{P^z_L, e_{x_0}}_{\car_L} = \inner{P^z_L, e_{x_0}}_{\cl_1 \car_L} = -P^z_L(x_0)
				$$
				for all $z \in \pa \car_L$. 
			\end{prop}
			\begin{proof}
				For $x \in \cl_1 \car_L$, let $K_x$ be an integer such that $x \in \pa \car_{K_x}$, if it exists. Let $P^x_{K_x}$ denote the Poisson kernel of $\car_{K_x}$ at boundary point $x \in \pa \car_{K_x}$. Define $e = e_{x_0} \in \R^{\cl_1 \car_L}$ by
				$$
				e(x) = -\epsilon(K_x) P^x_{K_x}(x_0)
				$$
				for some choice of values $\epsilon(k)$ such that $\sum_{k=0}^L \epsilon(k) = 1$ and $\epsilon(k) = 0$ if $k \leq K_{x_0}$. We will also choose $\epsilon(L) = 0$. If $x$ is such that no $K_x$ exists, define $e(x) = 0$.
				
				We claim that $\inner{P^z_L, e}_{\car_L} = \inner{P^z_L, e}_{\cl_1 \car_L}  = -P^z_L(x_0)$. To see this, let $X$ be a simple random walk on $\Z^d$ with law $\mathbb{P}^x$ started from $x$, and for $U \subset \Z^d$ let $\tau_U$ denote the first time $t$ that $X_t \not \in U$. We have
				\begin{align*}
					\inner{e, P^z_L}_{\cl_1 \car_L} &= -\sum_{k=K_{x_0}+1}^L \epsilon(k) \sum_{x \in \cl_1 \car_L : K_x = k} P^x_{k}(x_0) P^z_L(x) \\ &= -\sum_{k=K_{x_0}+1}^L \epsilon(k) \sum_{x \in \pa \car_k} \mathbb{P}^{x_0}(X_{\tau_{\car_K}} = x)\mathbb{P}^{x}(X_{\tau_{\car_L}} = z) \\ &= -\sum_{k=K_{x_0}+1}^L \epsilon(k) \sum_{x \in \pa \car_k} \mathbb{P}^{x_0}(X_{\tau_{\car_K}} = x, X_{\tau_{\car_L}} = z) \\
					&= -\sum_{k=0}^L \epsilon(k) P^z_L(x_0) = -P^z_L(x_0),
				\end{align*}
				where we have used the interpretation of the Poisson kernel as the exit location distribution of a simple random walk and the strong Markov property.
				
				We now choose $\epsilon$ to try to nearly minimize $\| e \|_{L^2(\cl_1 \car_L)}$. Using basic bounds on the Poisson kernel from \pref{greens_theo}, we have
				\begin{equation} \begin{aligned}
						\| e \|_{L^2(\cl_1 \car_L)}^2 &= \sum_{k=K_{x_0}  +1}^L |\epsilon(k)|^2 \sum_{x \in \cl_1 \car_L : K_x = K} |P^x_{k}(x_0)|^2 \\ &\leq C\sum_{k=K_{x_0}  +1}^L |\epsilon(k)|^2 \dist(x_0, \pa \car_k)^{-d+1} 
						\leq C \sum_{k = 1}^{ \rho_{x_0}} |\epsilon(k + K_{x_0})|^2 k^{-d+1}.
				\end{aligned} \end{equation}
				Note that
				$$
				\lambda := \sum_{k = 1}^{ \rho_{x_0}} k^{d/2 - 1/2} \geq C^{-1} \rho_{x_0}^{d/2 + 1/2},
				$$
				and so if we choose $\epsilon(k) = c_0 \lambda^{-1} \dist(x_0, \pa \car_k)^{d/2 - 1/2}$ for $K_{x_0} + 1 \leq k < L$ and some constant $c_0 \in [C^{-1}, C]$, we can ensure $\sum \epsilon(k) = 1$ and
				$$
				\| e \|_{L^2(\cl_1 \car_L)}^2 \leq C \lambda^{-2} \sum_{k=1}^{ \rho_{x_0}} 1 \leq C \rho_{x_0}^{-d}.
				$$
				This finishes the proof.
			\end{proof}
			
			\begin{remark}
				In the case that the domains $\car_L$ are replaced by smooth domains $\Omega_L$, we can carry out a similar construction and achieve the same estimates. Instead of defining $e_{x_0}$ with reference to discrete hypercubes $\car_{k}$ centered at $0$, we can center the hypercubes at $x_0$.
			\end{remark}
		\end{subsection}
		\begin{subsection}{Single-spin Law}
			The goal of this subsection is to prove the approximation $\mu^b_{x_0} \approx \nu^{b(x_0)}$ for $x_0$ not close to $\pa \car_L$. We do this through a change of variables and entropy argument. 
			
			Recall the definition of $\mathrm{J}(b(x_0))$ from \eref{nub_def}. For $f : \R \to \R$ measurable and bounded, we have
			$$
			\Esp_{\mu^b_{x_0}} f(\eta(x_0)) = \frac{\mathrm{J}(b(x_0))}{\K(\car_L, b)} \Esp_{\nu^{b(x_0)}} \left[ \K(\eta(x_0), x_0; b) f(\eta(x_0)) \right]
			$$
			where $\K(\car_L, b)$ is the partition function of $\Q^b$ and
			\begin{align*}
				\lefteqn{ \K(\xi, x_0; b) } \quad & \\ 
				& = \int_{\R^{\car_L}} \exp\(-\sum_{x \ne x_0} \( V(\eta(x)) - b(x)\eta(x) \) - \sum_{z \in \pa \car_L} V(\alpha^z(\eta)) \) \delta_{\xi}(d\eta(x_0))\prod_{x \ne x_0} d\eta(x).
			\end{align*}
			The measure $\delta_{\xi}(d\eta(x_0))$ is a Dirac delta enforcing $\eta(x_0) = \xi$. 
			
			The main task is to prove that $\K(\xi, x_0; b)$ depends only weakly on $\xi$. By the change of variables $\eta \mapsto \eta - \epsilon (\1_{x_0} + e_{x_0})$,
			\begin{align*}
				\lefteqn{ \K(\xi + \epsilon, x_0; b) } \quad & \\ &= \int_{\R^{\car_L}} e^{-\sum_{x \ne x_0} \( V((\eta + \epsilon e_{x_0})(x)) - b(x)(\eta + \epsilon e_{x_0})(x) \) - \sum_{z \in \pa \car_L} V(\alpha^z(\eta)) } d\delta_{\xi}(\eta(x_0))\prod_{x \ne x_0} d\eta(x).
			\end{align*}
			Thus
			\begin{equation} \label{e.dlogK}
				\frac{1}{\K(\xi,x_0;b)}\frac{d}{d\xi} \K(\xi,x_0; b) = -\Esp_{\Q^b} \left[ \sum_{x \ne x_0} V'(\eta(x)) e_{x_0}(x) - b(x)e_{x_0}(x) \ \bigg | \ \eta(x_0) = \xi \right].
			\end{equation}
			We use this computation in the next theorem, which shows that, far from the boundary, the one-spin marginal law does not feel strongly the effect of the boundary interaction.
			
			\begin{theo}
				The distance between $\mu^b_{x_0}$ and $\nu^{b(x_0)}$ in total variation norm is bounded by $C ( 1 + \| b \|_{L^2(\car_L)}) \rho_{x_0}^{-d/2}$. In particular, for any bounded and measurable $f$, we have
				\begin{equation} \label{e.TVdiff}
					| \Esp_{\mu^b_{x_0}} f - \Esp_{\nu^{b(x_0)}} f | \leq C \| f \|_{L^\infty(\R)} ( 1 + \| b \|_{L^2(\car_L)}) \rho_{x_0}^{-d/2}.
				\end{equation}
			\end{theo}
			\begin{proof}
				Without loss of generality, assume $x_0$ is not adjacent to $\pa \car_L$. Define 
				$$
				g(\xi) := \frac{\mu^b_{x_0}(d\xi)}{\nu^{b(x_0)}(d\xi)} = \frac{\mathrm{J}(b(x_0)) \K(\xi,x_0; b)}{\K(\car_L, b)}.
				$$ 
				Our strategy is to bound the relative entropy $H(\mu^{b}_{x_0} | \nu^{b(x_0)}) = \Esp_{\nu^{b(x_0)}} g(\xi) \log g(\xi)$. Pinsker's inequality will then allow us to bound the total variation distance.
				
				We apply the log-Sobolev inequality, as stated in \pref{poincarebe}, to see
				$$
				H(\mu^b_{x_0} | \nu^{b(x_0)}) := \Esp_{\nu^{b(x_0)}} g(\xi) \log g(\xi) \leq C\Esp_{\nu^{b(x_0)}} \left[\left( \frac{\pa_\xi g(\xi)}{\sqrt{g(\xi)}} \right)^2 \right].
				$$
				Let $\mathcal{I}(\xi) = (\K(\xi, x_0; b))^{-1} \frac{d}{d\xi} K(\xi, x_0; b)$. Then the RHS above can be written as
				\begin{equation}\begin{aligned}
						\Esp_{\nu^{b(x_0)}}  \left[ \left( \frac{\pa_\xi g(\xi)}{\sqrt{g(\xi)}} \right)^2 \right] &=\Esp_{\mu^{b}_{x_0}} \left[ \left( \frac{\pa_\xi g(\xi)}{{g(\xi)}} \right)^2 \right] =\Esp_{\mu^{b}_{x_0}} \left[ | \mathcal{I}(\xi) |^2 \right] \\ &\leq \Esp_{\Q^b} \left[\bigg | \sum_{x \ne x_0} V'(\eta(x)) e_{x_0}(x) - b(x) e_{x_0}(x) \bigg |^2 \right].
				\end{aligned}\end{equation}
				In the last inequality above, we used the expression \eref{dlogK} and Jensen's inequality. Continuing, we split the RHS above and apply Cauchy-Schwarz to get
				\begin{equation} \label{e.Hbd2}
					H(\mu^b_{x_0} | \nu^{b(x_0)}) \leq 2\Esp_{\Q^b} \left[\bigg| \sum_{x \ne x_0} V'(\eta(x)) e_{x_0}(x) \bigg|^2\right] + 2\| b \|_{L^2(\car_L)}^2 \| e_{x_0} \|_{L^2(\car_L)}^2.
				\end{equation}
				To bound the first term on the RHS, we first bound the mean 
				$$
				m(b) = \Esp_{\Q^b} \sum_{x \ne x_0} V'(\eta(x)) e_{x_0}(x).
				$$
				For any $r \in (0,1)$, we compute
				\begin{align}
					\frac{d}{dr} m(rb) &= \Cov_{\Q^{rb}} \left[ \sum_{x \in \car_L} b(x) \eta(x), \sum_{x \ne x_0} V'(\eta(x)) e_{x_0}(x) \right] \\ &\notag \leq  \(\Var_{\Q^{rb}} \left[ \sum_{x \in \car_L} b(x) \eta(x) \right] \)^{1/2} \( \Var_{\Q^{rb}}\left[ \sum_{x \ne x_0} V'(\eta(x)) e_{x_0}(x) \right] \)^{1/2}.
				\end{align}
				Now the Poincar\'e inequality implies $\Var_{\Q^{rb}} \left[  \sum_{x \in \car_L} b(x) \eta(x) \right] \leq C\| b \|_{L^2(\car_L)}^2$ and
				\begin{equation} \label{e.varbdry}
					\Var_{Q^{rb}}\left[\sum_{x \ne x_0} V'(\eta(x)) e_{x_0}(x) \right]\leq C \sum_{x \ne x_0} \Esp_{\Q^{rb}} \left[ |V''(\eta(x)) e_{x_0}(x)|^2 \right]\leq C \| e_{x_0} \|_{L^2(\car_L)}^2.
				\end{equation}
				Together with $m(0) = 0$ by symmetry of $\Q^0$ and $V$, we conclude 
				$$
				m(b) \leq C \| b \|_{L^2(\car_L)} \| e_{x_0} \|_{L^2(\car_L)}.$$
				The estimate \eref{varbdry} shows
				\begin{equation} \begin{aligned} \label{e.bdryint}
						\Esp_{\Q^b} \left[\bigg| \sum_{x \ne x_0} V'(\eta(x)) e_{x_0}(x) \bigg|^2\right] &= \Var_{\Q^b}\left[ \sum_{x \ne x_0} V'(\eta(x)) e_{x_0}(x) \right] + | m(b) |^2 \\ &\leq  C(1 + \| b \|_{L^2(\car_L)}^2) \| e_{x_0} \|_{L^2(\car_L)}^2.
				\end{aligned} \end{equation}
				We now apply this estimate to \eref{Hbd2} and use $\| e_{x_0} \|_{L^2(\car_L)}^2 \leq C \rho_{x_0}^{-d}$ to conclude 
				$$
				H(\mu^b_{x_0} | \nu^{b(x_0)}) \leq C (1 + \| b \|_{L^2(\car_L)}^2) \rho_{x_0}^{-d},
				$$
				and the total variation distance between $\mu^b_{x_0}$ and $\nu^{b(x_0)}$ is less than a constant times the square root of the relative entropy by Pinsker's inequality.
			\end{proof}
		\end{subsection}
		
		\begin{subsection}{Minimization of $\mathcal{F}_a(\cdot;\Q^b)$}
			Recall that $\ol{u}_a \in H^1(\Q^b; \R^{\car_L})$ is the minimizer of the functional $\mathcal{F}_{a}(\cdot; \Q^b)$ defined in \eref{Fadef}. By linearity of the map $a \mapsto \ol{u}_a$, it will be sufficient to understand $\ol{u}_{\1_x}$. Moreover, we have $\ol{u}_{\1_x}(y,\eta) = 0$ almost surely if $y \ne x$ by \pref{energy_bds}, so we will define $\ol{u}_x \in H^1(\Q^b; \R)$ by $\ol{u}_x(\eta) = \ol{u}_{\1_x}(x,\eta)$. By an abuse of notation, we will say $\ol{u}_x$ is the minimizer of $\mathcal{F}_{\1_x}(\cdot;\Q^b)$. Let $U^{b(x)} \in H^1(\nu^{b(x)};\R)$ be the minimizer of $\mathcal{G}_1(\cdot; \nu^{b(x)})$, which was defined in \eref{G0def}. We note that $\ol{u}_x$ depends on $b$, but we will mostly omit the dependence. However, all constants $C$ are independent of $b$. 
			
			The goal of this subsection is to prove $\ol{u}_x \approx U^{b(x)}$ in a certain sense, when $x$ is not close to $\pa \car_L$.  Recalling \pref{energy_bds}, we see that $\ol{u}_x$ solves
			$$
			-\Delta_\eta \ol{u}_x  + \sum_{y \in \car_L} \( V'(\eta(y)) - \sum_{z \in \pa \car_L} V'(\alpha^z(\eta)) P^z(y) - b(y) \) \pa_{\eta(y)} \ol{u}_x + V''(\eta(x)) \ol{u}_x = 1.
			$$
			The equation for $U^{b(x)}$ (as a function of $\eta(x)$) is
			$$
			-\Delta_{\eta(x)} U^{b(x)} + \( V'(\eta(x)) - b(x)\) \pa_{\eta(x)} U^{b(x)} + V''(\eta(x)) U^{b(x)} = 1.
			$$
			We define $w_x(\eta) := \ol{u}_x(\eta) - U^{b(x)}(\eta(x))$ and write the equation that $w_x$ solves:
			$$
			-\mathcal{L}_{\Q^b} w_x + V''(\eta(x)) w_x = - \sum_{z \in \pa \car_L} V'(\alpha^z(\eta)) P^z(x) \pa_{\eta(x)} U^{b(x)}.
			$$
			Integrating against $w_x \Q^b$ shows
			\begin{equation} \label{e.F0wx}
				2\mathcal{F}_0(w_x; \Q^b) = -\Esp_{\Q^b} \left[ w_x \sum_{z \in \pa \car_L} V'(\alpha^z(\eta)) P^z(x) \pa_{\eta(x)} U^{b(x)} \right].
			\end{equation}
			Our goal is to show that the RHS is smaller than $\mathcal{F}_0(w_x; \Q^b)$ plus some small constant. We cannot use H\"older's inequality directly on \eref{F0wx}, because in $d=2,3$ the summability properties of $P^z(x)$ are too weak (we can however get a strong enough result in $d \geq 4$). Instead, we must use the special profile $e_x \in \R^{\car_L}$ constructed in \pref{berg_con}. 
			
			\begin{lem}
				We have
				\begin{equation} \label{e.Fineq}
					\mathcal{F}_0(\ol{u}_{x_0} - U^{b(x_0)}; \Q^b) \leq C(1 + \|b\|_{L^2(\car_L)}^2) \rho_{x_0}^{-d}
				\end{equation}
				where $U^{b(x)}$ is a considered as a function of $\eta(x)$. Consequently
				\begin{equation} \label{e.sspde2}
					|\Esp_{\Q^b} \ol{u}_{x_0} - \Esp_{\nu^{b(x_0)}} U^{b(x_0)} | \leq C(1 + \|b\|_{L^2(\car_L)}) \rho_{x_0}^{-d/2}.
				\end{equation}
			\end{lem}
			\begin{proof}
				We will write $\Q^b(\eta)$ for the density of $\Q^b$ at $\eta$. Note that
				\begin{equation} \label{e.bdryQdecomp}
					\( \sum_{z \in \pa \car_L} V'(\alpha^z(\eta)) P^z(x_0) \) \Q^b(\eta) = -\pa_{\eta(x_0)} \Q^b(\eta) - V'(\eta(x_0)) \Q^b(\eta) + b(x_0) \Q^b(\eta).
				\end{equation}
				We break \eref{F0wx} into pieces
				\begin{equation} \begin{aligned}
						&\mathcal{I}_1 = \Esp_{\Q^b} \left[ w_{x_0} \pa_{\eta(x_0)} U^{b(x_0)} \frac{\pa_{\eta(x_0)} \Q^b(\eta)}{\Q^b(\eta)} \right] \\
						&\mathcal{I}_2 = \Esp_{\Q^b} \left[ w_{x_0} \pa_{\eta(x_0)} U^{b(x_0)} \(V'(\eta(x_0)) - b(x_0)\) \right]
				\end{aligned} \end{equation}
				so that $\mathcal{F}_0(w_{x_0}; \Q^b) = \mathcal{I}_1 + \mathcal{I}_2$.
				After performing a change of variables $\eta' = \eta + \epsilon e_{x_0}$, we have
				\begin{align}
					\lefteqn{ \Esp_{\Q^b} \left[ w_{x_0} \pa_{\eta(x_0)} U^{b(x_0)} \frac{\Q^b(\eta + \epsilon \1_{x_0})}{\Q^b(\eta)} \right] } \quad & \\ & \notag =\int_{\R^{\car_L}} w_{x_0}(\eta - \epsilon e_{x_0}) \pa_{\eta(x_0)} U^{b(x_0)}(\eta(x_0)) Q^b(\eta + \epsilon \1_{x_0} - \epsilon e_{x_0}) \prod_{x \in \car_L} d\eta(x)
				\end{align}
				for any $\epsilon \in \R$, and forming a difference quotient as $\epsilon \to 0$ shows that
				\begin{equation}
					\begin{aligned} \label{e.I1decomp}
						\mathcal{I}_1 &= -\Esp_{\Q^b} \left[ \(\sum_{x \ne x_0} e_{x_0}(x) \pa_{\eta(x)} w_{x_0} \) \pa_{\eta(x_0)} U^{b(x_0)} \right] \\ &\quad + \Esp_{\Q^b} \left[ w_{x_0} \pa_{\eta(x_0)} U^{b(x_0)} \( \sum_{x \ne x_0} (V'(\eta(x)) - b(x)) e_{x_0}(x) \) \right] \\ &\quad - \Esp_{\Q^b} \left[ w_{x_0} \pa_{\eta(x_0)} U^{b(x_0)} (V'(\eta(x_0)) - b(x_0)) \right].
					\end{aligned}
				\end{equation}
				The last line of \eref{I1decomp} cancels with $\mathcal{I}_2$, so we just need to bound the first two lines. On the first term use the $L^\infty$ bound in \pref{Genergy} on $\pa_{\eta(x_0)} U^{b(x_0)}$ to see
				\begin{align}
					\left| \Esp_{\Q^b} \left[ \(\sum_{x \ne x_0} e_{x_0}(x) \pa_{\eta(x)} w_{x_0} \) \pa_{\eta(x_0)} U^{b(x_0)} \right] \right| &\leq C \sum_{x \ne x_0} |e_{x_0}(x)| \Esp_{Q^b} {| \pa_{\eta(x)} w_{x_0} |} \\
					& \notag \leq \frac{1}{8} \sum_{x \ne x_0} \Esp_{\Q^b} |\pa_{\eta(x)} w_{x_0}|^2 + 2C \| e_{x_0} \|_{L^2(\car_L)}^2
				\end{align}
				Similarly on the second term we have
				\begin{multline}
					\left | \Esp_{\Q^b} \left[ w_{x_0} \pa_{\eta(x_0)} U^{b(x_0)} \( \sum_{x \ne x_0} (V'(\eta(x)) - b(x)) e_{x_0}(x) \) \right]  \right| \\ \leq \frac{1}{8 \inf V''} \Esp_{\Q^b} w_{x_0}^2 + C\Esp_{\Q^b} \left[ \( \sum_{x \ne x_0} V'(\eta(x)) e_{x_0}(x) \)^2 \right].
				\end{multline}
				Putting it all together into \eref{F0wx}, we get
				$$
				2\mathcal{F}_0(w_{x_0}; \Q^b) \leq \frac{1}{2} \mathcal{F}_0(w_{x_0}; \Q^b) + C \| e_{x_0} \|_{L^2(\car_L)}^2 + C\Esp_{\Q^b} \left[\( \sum_{x \ne x_0} V'(\eta(x)) e_{x_0}(x) \)^2 \right].
				$$
				The last term on the RHS was bounded by $C(1 + \| b \|_{L^2(\car_L)}^2) \| e_{x_0} \|_{L^2(\car_L)}^2$ in \eref{bdryint}, and so \eref{Fineq} is proved. For \eref{sspde2}, note that
				$$
				|\Esp_{\Q^b} \ol{u}_{x_0} - \Esp_{\mu^b_{x_0}} U^{b(x_0)} | = |\Esp_{\Q^b} w_{x_0}| \leq C\(\mathcal{F}_0(w_{x_0}; \Q^b)\)^{1/2},
				$$
				and, by \eref{TVdiff} and boundedness of $U^{b(x_0)}$, we have
				$$
				|\Esp_{\mu^b_{x_0}} U^{b(x_0)}  - \Esp_{\nu^{b(x_0)}} U^{b(x_0)}| \leq C(1 + \| b \|_{L^2(\car_L)}) \rho_{x_0}^{-d/2}.
				$$
				The triangle inequality finishes the proof.
			\end{proof}
			
			Since the map $a \mapsto \ol{u}_a$, where $\ol{u}_a$ minimizes $\mathcal{F}_a(\cdot;\Q^b)$, is linear, the lemma can be used to show $\ol{u}_a(x, \eta)$ and $U^{b(x)}(\eta) a(x)$ are approximately equal so long as $a$ does not have its mass concentrate near $\pa \car_L$.
			
			To conclude the section, we consider the dependence of $U^{b(x_0)}$ on $b$.
			\begin{lem}
				Let $U^{\beta}$ denote the minimizer of $\mathcal{G}_1(\cdot; \nu^{\beta})$, and let $U^0$ denote the minimizer of $\mathcal{G}_1(\cdot; \nu^0)$. We have
				\begin{equation} \label{e.bpert}
					| \Esp_{\nu^{\beta}} [ U^\beta ] - \Esp_{\nu^0} [ U^0]  | \leq C \beta
				\end{equation}
				for some constant $C$.
			\end{lem}
			\begin{proof}
				First, the difference $w := U^\beta - U^0$ satisfies the equation
				$$
				-\Delta w + V'(\xi) \pa w + V''(\xi) w = \beta \pa U^\beta
				$$
				over $\R$. Integrating against $w \nu^0$, and using the bound on $\pa U^\beta$ from \pref{Genergy}, we obtain
				$$
				\mathcal{G}_0(w; \nu^0) \leq C |\beta| \Esp_{\nu^0}(|w|) \leq C |\beta| \sqrt{\mathcal{G}_0(w; \nu^0)} \leq C\beta^2,
				$$
				and so
				$$
				\left | \Esp_{\nu^0} U^\beta  - \Esp_{\nu^0} U^0 \right |^2 \leq \Esp_{\nu^0} |U^\beta - U^0|^2 \leq C\mathcal{G}_0(w; \nu^0) \leq C \beta^2.
				$$
				By interpolation and the Poincar\'e inequality, we have
				$$
				|\Esp_{\nu^\beta} U^\beta - \Esp_{\nu^0} U^\beta | \leq \int_0^1 \left|\Cov_{\nu^{r\beta}} \left[ \beta \xi, U^\beta(\xi) \right] \right| dr \leq C|\beta| \( \Var_{\nu^{r\beta}} U^\beta \)^{1/2} \leq C |\beta|.
				$$
				We conclude with the triangle inequality.
			\end{proof}
		\end{subsection}
		
	\end{section}

	\begin{section}{Minimization of the Helffer-Sj\"ostrand Energy} \label{s.Vpert_sec}
		In this section, we will take general $a \in \R^{\car_L}$, but we have particular interest in the case where $a$ is almost orthogonal to the family $\{ P^z \}_{z \in \pa \car_L}$ in the $L^2(\car_L)$ inner product. We consider $(\eta(x))_{x \in \car_L}$ distributed according to the measure $\Q^b$. For any $x \in \car_L$, we let $\mu^b_x$ be the marginal distribution of $\eta(x)$ under $\Q^b$, and let $\nu^{b(x)}$ be the probability measure on $\R$ proportional to $\exp(-V(\xi) + b(x) \xi)d\xi$.
		
		We seek a precise understanding of $\mathcal{E}_a(\cdot;\Q^b)$ and its minimizer. In the case of the Gaussian membrane model and $\inner{a,P^z}_{\car_L}= 0$, the minimizer is simply $u_a(x,\eta) = c a(x)$ for some constant $c$. Critically, the fact that $a$ is "boundary orthogonal", i.e.\  orthogonal to the $\{ P^z \}$, means that the boundary terms within $\mathcal{E}_a$ are irrelevant for the Gaussian energy.
		
		In the non-Gaussian case, even if $\inner{a,P^z}_{\car_L} = 0$, the boundary terms are still active due to the random background measure $\Q^b$ and random coefficients involving $V''$ within $\mathcal{E}_a$. We pursue the idea, mentioned at the end of \sref{intro_sec}, that there should be a homogenization or law of large numbers effect due to the random environment $\eta$ in the boundary terms within $\mathcal{E}_a$. It will allow us to show the minimizer of $\mathcal{E}_a$ is close to that of $\mathcal{F}_a$, which was understood in \sref{SSlaw_sec}.
		
		The question of dependence of $u_a = \mathrm{HS}^b a$ on $b$ must also be addressed. We find that if $|a|^2 |b|^2$ is small in $L^1(\car_L)$, the dependence is also small. This condition will hold for the limits involving the rescaled field $\ol{\phi}$, but not for the infinite volume limit in $d \geq 5$. For the infinite volume limit \tref{inf_limit_thm}, we will consider $a$ and $b$ with mass concentrated on a single point in \pref{pt_est_theo}.
		
		Fix $\epsilon \in (0,1)$ and $\ell \in \Z$ such $L^{1-\epsilon} \leq \ell \leq 2L^{1-\epsilon}$, and let $\Lambda = \car_L \setminus \car_{L - \ell}$ be a boundary layer of thickness $\ell$. Define $w \in H^1(\Q^b; \R^{\car_L})$ by
		\begin{equation} \label{e.w_min_def}
			w(x,\eta) = \begin{cases}
				U^{b(x)}(\eta(x)) a(x), & \quad x \in \car_{L - \ell}, \\
				(\Esp_{\nu_0} U^0)a(x), &\quad x \in \Lambda.
			\end{cases}
		\end{equation}
		Here $U^{b(x)}$ is the minimizer of $\mathcal{G}_1(\cdot;\nu^{b(x)})$. We claim that $w$ is an approximate minimizer of $\mathcal{E}_a$, and therefore a good approximation to the true minimizer $u_a$ by coercivity of $\mathcal{E}_0$. We address the boundary energy of $w$ first.
		\begin{lem} \label{l.bdry_est_lem}
			With $w$ as in \eref{w_min_def}, there is a constant $C > 0$ such that
			\begin{align} \label{e.bdry_est}
				\lefteqn{\Esp_{\Q^b}  \sum_{z \in \pa \car_L} V''(\alpha^z(\eta)) |\inner{P^z, w}_{\car_L}|^2} \quad & \\ \notag &\leq C(1 + \| b \|_{L^2(\car_L)}^2) \| a \|_{L^2(\car_L)}^2 L^{\epsilon(3d-2)}L^{-(d-1)} + C \cdot B_{L}(a),
			\end{align}
			where
			\begin{equation} \label{e.BL_def}
				B_{L}(a) :=  \sum_{z \in \pa \car_L} |\inner{P^z,a}_{\car_L}|^2.
			\end{equation}
			is the "boundary energy" of $a$.
		\end{lem}
		\begin{proof}
			For $z \in \pa \car_L$, we decompose $\Esp_{\Q^b} \inner{P^z, w}_{\car_L}$ into two pieces:
			$$
			\Esp_{\Q^b} \inner{P^z, w}_{\car_L} = \Esp_{\nu^0} U^0 \inner{P^z, a}_{\car_L} +  \sum_{x \in \car_{L - \ell}} P^z(x) (\Esp_{\mu^b_x} U^{b(x)} - \Esp_{\nu^0} U^0) a(x).
			$$
			Applying \eref{TVdiff} for $x \in \car_{L - \ell}$, we see
			$$
			| \Esp_{\mu^b_x} U^{b(x)} - \Esp_{\nu^{b(x)}} U^{b(x)} |\leq C(1 + \| b \|_{L^2(\car_L)}) \ell^{-d/2}
			$$
			and from \eref{bpert} we have
			$$
			|\Esp_{\nu^{b(x)}} U^{b(x)} - \Esp_{\nu^0} U^0| \leq C |b(x)|.
			$$
			Thus, bounding $|P^z(x)| \leq C\ell^{-(d-1)}$ (by \eref{poisson}) and using Cauchy-Schwarz, we have
			\begin{align}
				\lefteqn{| \Esp_{\Q^b} \inner{P^z, w}_{\car_L} - \Esp_{\nu^0} U^0 \inner{P^z, a}_{\car_L} |} \quad &\\ \notag &\leq C\sum_{x \in \car_{L - \ell}} \ell^{-d+1} (|b(x)| + (1+ \|b\|_{L^2(\car_L)}) \ell^{-d/2}) |a(x)| \\ \notag
				&\leq C\ell^{-(d-1)} \| b \|_{L^2(\car_L)} \|a \|_{L^2(\car_L)}  + C(1 + \|b\|_{L^2(\car_L)})\ell^{-\frac32 d + 1} \| a \|_{L^1(\car_L)} \\ \notag
				&\leq C\ell^{-(d-1)} \| b \|_{L^2(\car_L)} \|a \|_{L^2(\car_L)} + C(1+\|b\|_{L^2(\car_L)}) \ell^{-(d-1)} L^{\epsilon d/2} \|a\|_{L^2(\car_L)} \\ \notag
				&\leq CL^{\epsilon (3d/2 - 1)} L^{-(d-1)}(1 + \| b \|_{L^2(\car_L)})\|a \|_{L^2(\car_L)} =: \mathrm{Error}_1.
			\end{align}
			We conclude $|\Esp_{\Q^b} \inner{P^z, w}|^2 \leq 2(\mathrm{Error}_1)^2 + C|\inner{P^z,a}_{\car_L}|^2$.
			
			Next, we look at the variances of $\inner{P^z,w}_{\car_L}$ under $\Q^b$. Applying the Poincar\'e inequality for $\Q^b$ gives
			\begin{align}
				\sum_{z \in \pa \car_L} \Var_{\Q^b} \inner{P^z, w}_{\car_L} &\leq C \Esp_{\Q^b} \sum_{x \in \car_{L - \ell}} \( \sum_{z \in \pa \car_L} |P^z(x)|^2 \) |a(x)|^2|\pa_{\eta(x)} U^{b(x)}(\eta(x))|^2\\  & \notag \leq C \ell^{-d+1} \| a \|_{L^2(\car_L)}^2 \sup_{x \in \car_{L - \ell}} \Esp_{\Q^b} | \pa_{\eta(x)} U^{b(x)}(\eta(x)) |^2 \\ &\notag \leq C L^{\epsilon(d-1)}L^{-(d-1)} \| a \|_{L^2(\car_L)}^2.
			\end{align}
			Here we used the estimate
			$$
			\sum_{z \in \pa \car_L}|P^z(x)|^2 \leq \sup_{z \in \pa \car_L} P^z(x) \sum_{z \in \pa \car_L} P^z(x) = \sup_{z \in \pa \car_L} P^z(x) \leq C \ell^{-d+1}.
			$$
			Thus the boundary term in the energy is
			\begin{align}
				\lefteqn{\frac12\Esp_{\Q^b}  \sum_{z \in \pa \car_L} V''(\alpha^z(\eta)) |\inner{P^z, w}_{\car_L}|^2} \quad & \\ \notag &\leq C \sum_z \( \Var_{\Q^b} \inner{P^z, w}_{\car_L} + | \Esp_{\Q^b} \inner{P^z, w}_{\car_L}|^2 \) \\ \notag &\leq C B_L(a) + CL^{d-1}(\text{Error}_1)^2 + C L^{\epsilon(d-1)}L^{-(d-1)} \| a \|_{L^2(\car_L)}^2 \\ \notag &\leq C B_L(a) + C(1 + \| b \|_{L^2(\car_L)}^2) \| a \|_{L^2(\car_L)}^2 L^{\epsilon(3d-2)}L^{-(d-1)}.
			\end{align}
			This completes the proof.
		\end{proof}
		
		Next, we show that $w$ nearly minimizes $\mathcal{F}_a$. 
		\begin{lem} \label{l.Fapproxlem}
			We have
			\begin{equation}
				\mathcal{F}_a(w; \Q^b) \leq \mathcal{F}_a(\ol{u}_a; \Q^b) + \mathrm{Error}
			\end{equation}
			where
			$$
			|\mathrm{Error}| \leq C \| a \|_{L^2(\Lambda)}^2 + C(1 + \| b \|_{L^2(\car_L)}^2)\| a \|_{L^2(\car_L)}^2L^{\epsilon d} L^{-d}.
			$$
		\end{lem}
		\begin{proof}
			Note that $\ol{u}_a(x,\eta) = a(x)\ol{u}_x(\eta)$, where $\ol{u}_x \1_x$ is the minimizer of $\mathcal{F}_{\1_{x}}(\cdot;\Q^b)$. We have
			\begin{equation} \begin{aligned}
					\mathcal{F}_a(w; \Q^b) - \mathcal{F}_a(\ol{u}_a; \Q^b) &= \mathcal{F}_0(w - \ol{u}_a; \Q^b) \\ &= \sum_{x \in \Lambda} |a(x)|^2 \mathcal{F}_0((\Esp_{\nu^0} U^0) \1_{x} - \ol{u}_x \1_x; \Q^b)
					\\ &\quad +\sum_{x \in \car_{L - \ell}} |a(x)|^2\mathcal{F}_0(U^{b(x)} \1_x - \ol{u}_x \1_x; \Q^b).
			\end{aligned} \end{equation}
			We now control the two sums on the RHS above, which we label $S_1$ and $S_2$, respectively. For the first, we have
			$$
			\mathcal{F}_0((\Esp_{\nu^0} U^0) \1_{x} - \ol{u}_x \1_x; \Q^b) \leq 2\mathcal{F}_0(\ol{u}_x \1_x; \Q^b) + 2\mathcal{F}_0((\Esp_{\nu^0} U^0) \1_{x}; \Q^b) \leq C,
			$$
			and so $S_1 \leq C\| a \|_{L^2(\Lambda)}^2$. The term $S_2$ is controlled by \eref{Fineq}. Specifically, we have
			$$
			\mathcal{F}_0(U^{b(x)} \1_x - \ol{u}_x \1_x; \Q^b) \leq C(1 + \| b \|_{L^2(\car_L)}^2) \ell^{-d}
			$$
			and inserting this bound into $S_2$ finishes the proof.
		\end{proof}
		
		The main result of the section follows.
		\begin{prop}
			Let $w$ be defined as in \eref{w_min_def}, and recall the definition of $B_L(a)$ in \eref{BL_def}. For $u_a = \mathrm{HS}^b a$, we have
			\begin{align} \label{e.E0_nearmin}
				\mathcal{E}_0(w - u_a; \Q^b) &\leq C \| a \|_{L^2(\Lambda)}^2 + C B_L(a) \\ \notag &\quad + C(1 + \| b \|_{L^2(\car_L)}^2)\| a \|_{L^2(\car_L)}^2 L^{(3d-2)\epsilon} L^{-(d-1)},
			\end{align}
			and
			\begin{align} \label{e.mean_scaling}
				\| \Esp_{\Q^b} u_a - (\Esp_{\nu^0} U^0)a \|_{L^2(\car_L)} &\leq C(1 + \| b \|_{L^2(\car_L)}) \| a \|_{L^2(\car_L)} L^{\epsilon(3d/2 - 1)}L^{-(d-1)/2} \\ \notag &+ C \| a \|_{L^2(\Lambda)} + C({B_L(a)})^{1/2} + C \( \sum_{x \in \car_L} |a(x)|^2 |b(x)|^2\)^{1/2}.
			\end{align}
			In the case that $\| a \|_{L^2(\Lambda)}^2 \leq CL^{-\epsilon}$, $\| a \|_{L^2(\car_L)} \leq C$, $\sup |b(x)| \leq CL^{-d/2}$ and $B_L(a) \leq CL^{-1}$, we get
			$$
			\| \Esp_{\Q^b} u_a - (\Esp_{\nu^0} U^0)a \|_{L^2(\car_L)} \leq CL^{\epsilon(3d/2 - 1)}L^{-(d-1)/2} + CL^{-\epsilon/2} + CL^{-d/2} + CL^{-1/2}.
			$$
			We can then optimize by choosing $\epsilon = \frac{d-1}{3d-1}$ to get a rate of $CL^{-\frac{d-1}{6d-2}}$ on the RHS.
		\end{prop}
		\begin{proof}
			Using quadratic response and $\inf \mathcal{F}_a \leq \inf \mathcal{E}_a$, we start from
			$$
			\mathcal{E}_0(w - u_a; \Q^b) = \mathcal{E}_a(w; \Q^b) - \mathcal{E}_a(u_a; \Q^b) \leq \mathcal{E}_a(w; \Q^b) - \mathcal{F}_a(\ol{u}_a; \Q^b).
			$$
			We bound the boundary term inside $\mathcal{E}_a(w; \Q^b)$ using \lref{bdry_est_lem}, and \lref{Fapproxlem} bounds the remainder:
			\begin{equation} \begin{aligned}
					\mathcal{E}_a(w;\Q^b) &\leq \mathcal{F}_a(w;\Q^b) + C(1 + \| b \|^2_{L^2(\car_L)}) \| a \|_{L^2(\car_L)}^2 L^{\epsilon(3d-2)} L^{-(d-1)} + C B_L(a)\\
					&\leq \mathcal{F}_a(\ol{u}_a; \Q^b) + C \| a \|_{L^2(\Lambda)}^2  \\ &\quad +C(1 + \| b \|^2_{L^2(\car_L)}) \| a \|_{L^2(\car_L)}^2 L^{\epsilon(3d-2)} L^{-(d-1)} + C B_L(a),
			\end{aligned} \end{equation}
			and so
			$$
			\mathcal{E}_0(w - u_a; \Q^b) \leq C \| a \|_{L^2(\Lambda)}^2 + C(1 + \| b \|_{L^2(\car_L)}^2) \| a \|_{L^2(\car_L)}^2 L^{\epsilon(3d-2)} L^{-(d-1)} + C B_L(a).
			$$
			Clearly,
			$$
			\sum_{x \in \car_{L}} |\Esp_{\Q^b} u_a(x) - \Esp_{\Q^b}w(x)|^2 \leq C\mathcal{E}_0(w - u_a; \Q^b),
			$$
			and by \eref{TVdiff}, which shows $\nu^{b(x)}$ approximates the marginal of $\Q^b$ at $\eta(x)$, we have
			$$
			\sum_{x \in \car_{L}} |\Esp_{\Q^b}w(x) - \Esp_{\nu^{b(x)}} w(x) |^2 \leq C(1 + \| b \|_{L^2(\car_L)}^2) \| a \|_{L^2(\car_L)}^2 \ell^{-d}.
			$$
			Finally, by \eref{bpert}
			$$
			\sum_{x \in \car_{L}} |\Esp_{\nu^{b(x)}}w(x) - a(x) \Esp_{\nu^0} U^0 |^2 \leq C\sum_{x \in \car_L} |a(x)|^2 |b(x)|^2.
			$$
			Combining the last three displays and \eref{E0_nearmin} with the triangle inequality proves \eref{mean_scaling}.
		\end{proof}
		
		The bound \eref{mean_scaling} is satisfactory when the perturbations $a, b$ have the following properties: (1) they are bounded in $L^2$, (2) $a$ does not concentrate on a thin boundary layer, and (3) $|a|^2|b|$ is small in $L^1$. Condition (3) does not hold in the infinite volume limit \tref{inf_limit_thm}, so we now provide a different construction that works in this case. 
		
		It will turn out to be sufficient to assume that $a$ concentrates almost all of its mass on a single point $x_0$. In this case, the minimizer $u_a$ of $\mathcal{E}_a$ is similarly concentrated on $x_0$, and there is a simple approximation for $\Esp_{\Q^b} u_a(x_0,\eta)$ given by the energy $\mathcal{G}_1$.
		
		\begin{prop} \label{p.pt_est_theo}
			Let $a \in \R^{\car_L}$ and let $u_a = \mathrm{HS}^b a$, and recall the definition of $B_L(a)$ from \eref{BL_def}. We have
			\begin{equation} \label{e.ptmass}
				\Esp_{\Q^{b}} \| u_a \|_{L^2(\car_L \setminus \{ x_0 \})}^2 \leq C \| a \|_{L^2(\car_L \setminus \{ x_0 \})}^2 + CB_L(a),
			\end{equation}
			and
			\begin{align} \label{e.meanthermo}
				\lefteqn{ | \Esp_{\Q^b}u_a(x_0) - a(x_0) \Esp_{\nu^{b(x_0)}} U^{b(x_0)} |^2 } \quad & \\ \notag & \leq C (1 + \| b \|_{L^2(\car_L)}^2) |a(x_0)|^2\rho_{x_0}^{-d} + C \| a \|_{L^2(\car_L \setminus \{ x_0\})}^2 + CB_L(a).
			\end{align}
		\end{prop}
		\begin{proof}
			Let $w(x,\eta) = \ol{u}_{\1_{x_0}}(x_0,\eta)a(x)$ for all $x \in \car_L$, where $\ol{u}_{\1_{x_0}}$ minimizes $\mathcal{F}_{\1_{x_0}}(\cdot;\Q^b)$. Observe that
			\begin{equation*}
				\Esp_{\Q^b} \sum_{z \in \pa \car_L} V''(\alpha^z(\eta)) |\inner{P^z, w}_{\car_L}|^2 =  \sum_{z \in \pa \car_L} \Esp_{\Q^b} \left[ V''(\alpha^z(\eta)) |\ol{u}_{\1_{x_0}}(x_0,\eta)|^2 \right] |\inner{P^z, a}_{\car_L}|^2,
			\end{equation*}
			and this is bounded by $C B_L(a)$ by \eref{upbd1}. We can apply the above to see
			$$
			\mathcal{E}_a(w; \Q^b) - \mathcal{E}_a(u_a; \Q^b) \leq \mathcal{E}_a(w; \Q^b) - \mathcal{F}_a(\ol{u}_a; \Q^b) \leq \mathcal{F}_a(w; \Q^b) - \mathcal{F}_a(\ol{u}_a; \Q^b) + CB_L(a).
			$$
			Next, we have
			$$
			\mathcal{F}_a(w; \Q^b) - \mathcal{F}_a(\ol{u}_a; \Q^b) = \mathcal{F}_0(w - \ol{u}_a; \Q^b) = \sum_{x \in \car_L} |a(x)|^2 \mathcal{F}_0(\ol{u}_{\1_{x_0}}(x_0,\cdot) \1_x(\cdot) - \ol{u}_{\1_x}; \Q^b).
			$$
			For $x \ne x_0$, we bound
			$$
			\mathcal{F}_0(\ol{u}_{\1_{x_0}}(x_0,\cdot) \1_x(\cdot) - \ol{u}_{\1_x}; \Q^b) \leq C
			$$
			and for $x = x_0$, we have
			$$
			\mathcal{F}_0(\ol{u}_{\1_{x_0}}(x_0,\cdot) \1_x(\cdot) - \ol{u}_{\1_x}; \Q^b) = 0.
			$$
			due to the fact $\ol{u}_{\1_x}(y,\eta) = 0$ a.s.\ for $y \ne x$. Thus
			$$
			\mathcal{F}_0(w - \ol{u}_a; \Q^b) \leq C \| a \|_{L^2(\car_L \setminus \{ x_0 \})}^2.
			$$
			It follows that
			\begin{align}
				\Esp_{\Q^{b}} \| u_a \|_{L^2(\car_L \setminus \{ x_0 \})}^2 &\leq C\Esp_{\Q^b} \| w \|_{L^2(\car_L \setminus \{ x_0 \})}^2 +  C\mathcal{E}_0(w - u_a; \Q^b) \\ \notag &\leq C \| a \|_{L^2(\car_L \setminus \{ x_0 \})}^2 + CB_L(a),
			\end{align}
			as desired.
			
			We now prove \eref{meanthermo}. We compute, by the triangle inequality and Jensen's inequality, that
			\begin{align*}
				| \Esp_{\Q^b}u_a(x_0,\eta) - a(x_0) \Esp_{\nu^{b(x_0)}} U^{b(x_0)} |^2 &\leq C\Esp_{\Q^b}\| u_a - w\|_{L^2(\car_L)}^2 \\ &\quad + C|a(x_0)|^2| \Esp_{\Q^b}\ol{u}_{x_0} - \Esp_{\nu^{b(x_0)}} U^{b(x_0)} |^2.
			\end{align*}
			The first term on the RHS above is dominated by $C\mathcal{E}_0(w - u_a; \Q^b)$, which we have already estimated. By \eref{sspde2}, we can bound the second term by
			$$
			| \Esp_{\Q^b}\ol{u}_{\1_{x_0}} - \Esp_{\nu^{b(x_0)}} U^{b(x_0)} |^2 \leq C(1 + \| b \|^2_{L^2(\car_L)}) \rho_{x_0}^{-d}.
			$$
			The proof is finished.
		\end{proof}
		
	\end{section}
	
	\begin{section}{Proofs of the Main Theorems} \label{s.concl_sec}
		In this section, we prove the main theorems listed in \sref{intro_sec}. The general process for all the results is as follows.
		\begin{enumerate}
			\item Figure out the relevant $a \in \R^{\car_L}$ at which we need to compute the cumulant generating function (c.g.f.) of $\Q$.
			\item Break $a$ down into a harmonic part $K_L a$ and an "almost boundary-orthogonal" part $K_L^\perp a$ with the Bergman projection $K_L$, which is the $L^2(\cl_1 \car_L)$ projection onto the space of functions which are discrete harmonic on $\car_L$.
			\item Estimate the sizes of $a$, $K_L a$, $K_L^\perp a$, especially in a boundary layer of $\car_L \cup \pa \car_L$. The estimates for $K_L a$ and $K_L^\perp a$ are major tasks in \sref{append}.
			\item Apply the estimates from \sref{harm_pert_sec} and \sref{Vpert_sec} to estimate the solutions of Helffer-Sj\"ostrand equations.
			\item Apply the Helffer-Sj\"ostrand representation of the c.g.f. given by \eref{cumulant_var} and \eref{var_form}.
			\item  For the limits of $\ol{\phi}$, this process applies equally well to the Gaussian membrane model, so we can compare the non-Gaussian to the Gaussian case and achieve our results.
		\end{enumerate}
		
		\begin{subsection}{Infinite Volume Limit} We first prove the infinite volume limit \tref{inf_limit_thm} in $d \geq 5$. Fix a $K \in \N$ to be held constant as $L \to \infty$, and let $a' \in \R^{\Z^d}$ be supported within $\car_K$. Let $a \in \R^{\car_L}$ be defined by $a(x) = -\sum_{y \in \car_K} \Gamma_L(x,y) a'(y)$, where $-\Gamma_L$ is the Dirichlet Green's function for $\Delta$ on $\car_L$. We allow constants denoted by $C$ to depend on $K$ and $\| a' \|_{L^2(\car_K)}$, but they remain uniform as $L \to \infty$.
			
			We let $K_L$ be the $L^2(\car_L \cup \pa \car_L)$ projection onto the linear span of $\{P^z\}_{z \in \pa \car_L}$ and $K_L^\perp = \text{Id} - K_L$. These operators are analyzed in \sref{append}.
			
			\begin{theo}
				For any $\epsilon > 0$ fixed, $d \geq 5$, $r \in [0,1]$, and $a$ as above, we have
				\begin{equation} \label{e.thermoform}
					\Esp_{\Q^{ra}} \inner{a, \mathrm{HS}^{ra} a}_{\car_L} = \sum_{x \in \car_L} |a(x)|^2 \Esp_{\nu^{ra(x)}} U^{ra(x)} + O(L^{-d/2 + 2 + \epsilon })
				\end{equation}
				as $L \to \infty$. Consequently, \tref{inf_limit_thm} holds, characterizing the infinite volume limit of the membrane model.
			\end{theo}
			\begin{proof}
				We will first prove \eref{thermoform} and then show how it implies \tref{inf_limit_thm}.
				
				We give some routine bounds on $a$. By \pref{greens_theo}, we have
				\begin{equation} \begin{aligned}
						\| a \|_{L^2(\car_L)}^2 &\leq C \sum_{x \in \car_L} \sum_{y \in \car_K} |a'(y)|^2 |\Gamma_L(x, y)|^2 \\ &\leq C \sum_{x \in \Z^d} \sup_{y \in \car_K}  |\Gamma(x,y)|^2 \leq C \sum_{x \in \Z^d} \frac{1}{(1 + |x|)^{2d-4}} \leq C.
				\end{aligned}\end{equation}
				Next we estimate $a$ on a boundary layer. Fix $\epsilon \in (0,1/4)$, and let $\ell \in \Z$ be between $L^{1-\epsilon}$ and $2L^{1-\epsilon}$. Let $\Lambda = \car_L \setminus \car_{L -\ell}$ be a boundary layer of thickness $\ell$. Then the same computation as before shows
				\begin{equation} \label{e.abdlayer}
					\| a \|_{L^2(\Lambda)}^2 \leq C \sum_{x \in \Lambda} \sum_{y \in \car_K} |a'(y)|^2 |\Gamma_L(x, y)|^2 \leq C \sum_{x \in \Lambda} \frac{1}{|x|^{2d-4}} \leq C L^{-d + 4 - \epsilon}.
				\end{equation}
				Similarly we bound the $L^1(\car_L)$ norm of $a$ by
				$$
				\| a \|_{L^1(\car_L)} \leq C \sup_{y \in \car_K} \sum_{x \in \car_L} \frac{1}{ (1 + |x - y|)^{d-2}}\leq C L^2.
				$$
				
				We now compute $\Esp_{\Q^{ra}} \inner{a, \mathrm{HS}^{ra} a}_{\car_L}$. Recall the boundary energy $B_L(a)$ defined in \eref{BL_def} and the special profile $e_x$ constructed in \pref{berg_con}. Recall also that $\rho_x$ denotes the distance from $x$ to $\pa \car_L$. 
				
				We use the decomposition
				\begin{equation} \label{e.decomp_infvol}
					\mathrm{HS}^{ra} a = \mathrm{HS}^{ra} (a \1_{\Lambda}) + \sum_{x \in \car_{L - \ell}} a(x) \mathrm{HS}^{ra} ({K}_L \1_{x}) + \sum_{x \in \car_{L - \ell}} a(x) \mathrm{HS}^{ra} ({K}_L^\perp \1_{x}),
				\end{equation}
				and accordingly define
				\begin{equation} \begin{aligned}
					T_1 &= \inner{a, \Esp_{\Q^{ra}} \mathrm{HS}^{ra}(a \1_\Lambda)} \\
					T_2 &= \sum_{x \in \car_{L - \ell}} a(x) \inner{a, \Esp_{\Q^{ra}}\mathrm{HS}^{ra}(K_L \1_x)} \\
					T_3 &= \sum_{x \in \car_{L - \ell}} a(x) \inner{a, \Esp_{\Q^{ra}} \mathrm{HS}^{ra}(K^\perp_L \1_x)}.
				\end{aligned} \end{equation}
				Only the last term $T_3$ should make a non-negligible contribution toward \eref{thermoform}. Indeed, the first term $T_1$ is a boundary layer term. For $T_2$, we will use the special profile $e_x$ to show that $K_L \1_x$ is very small in $L^2(\car_L)$ when $x$ is not close to $\pa \car_L$. The term $T_3$ will be simplified further with the approximation $$
				\Esp_{\Q^{ra}}\left[ \mathrm{HS}^{ra}(K_L^\perp \1_x)(y,\eta) \right] \approx a(x) \( \Esp_{\nu^{ra(x)}} U^{ra(x)} \) \1_{x}(y).
				$$
				
				We handle the first term by
				\begin{equation*}\begin{aligned} \label{e.T1bd}
						|T_1| \leq \| a \|_{L^2(\car_L)} \| \Esp_{\Q^{ra}} \mathrm{HS}^{ra}(a \1_{\Lambda}) \|_{L^2(\car_L)} \leq C \| a \|_{L^2(\car_L)} \| a \1_{\Lambda}\|_{L^2(\car_L)} \leq C L^{-d/2 + 2 -\epsilon/2},
				\end{aligned}\end{equation*}
				where we used that $\Esp_{\Q^{ra}} \mathrm{HS}^{ra}$ is bounded as an operator on $L^2(\car_L)$ uniformly in $L$. Indeed, by \eref{upbd1} and dominating the squared $L^2(\Q^{ra}; \R^{\car_L})$ norm by the energy $C\mathcal{E}_0(\cdot;\Q^{ra})$, we see that $\mathrm{HS}^{ra}$ is bounded from $L^2(\car_L)$ to $L^2(\Q^{ra}; \R^{\car_L})$. By Jensen's inequality, we have for any $b' \in \R^{\car_L}$ that
				$$
				\| \Esp_{\Q^{ra}} \mathrm{HS}^{ra} b' \|_{L^2(\car_L)}^2 \leq \Esp_{\Q^{ra}} \| \mathrm{HS}^{ra} b'(\cdot, \eta) \|_{L^2(\car_L)}^2 =: \| \mathrm{HS}^{ra} b' \|_{L^2(\Q^b ; \R^{\car_L})}^2,
				$$
				and so $\Esp_{\Q^{ra}} \mathrm{HS}^{ra}$ is bounded on $L^2(\car_L)$.
				
				For the second term $T_2$, observe that ${K}_L(\1_x + e_x) = 0$ since $\inner{P^z, \1_x + e_x}_{\cl_1(\car_L)} = 0$ for all $z \in \pa \car_L$ by construction. We have 
				\begin{equation} \label{e.berg_infvol}
					\| {K}_L \1_{x} \|_{L^2(\cl_1 \car_L)} = \| {K}_L e_x \|_{L^2(\cl_1 \car_L)} \leq C \rho_{x}^{-d/2},
				\end{equation}
				where we used that ${K}_L$ is a projection operator and $\| e_x \|_{L^2(\cl_1 \car_L)} \leq C \rho_x^{-d/2}$. Applying Cauchy-Schwarz and boundedness of $\Esp_{\Q^{ra}} \mathrm{HS}^{ra}$ shows
				\begin{equation} \begin{aligned}
						|T_2| &\leq C  \sum_{x \in \car_{L - \ell}}|a(x)| \| a \|_{L^2(\car_L)}  \rho_x^{-d/2} \\ &\leq C \|a \|_{L^1(\car_L)} \| a \|_{L^2(\car_L)} \ell^{-d/2} \leq C L^{\epsilon d/2} L^{-d/2 + 2}.
				\end{aligned} \end{equation}
				The term $T_2$ is thus negligible as $L \to \infty$ provided $\epsilon$ is chosen small.
				
				For $T_3$, we write $\inner{a, \Esp_{\Q^{ra}} \mathrm{HS}^{ra}{K}_L^\perp \1_{x}}_{\car_L}$ as $a(x) \Esp_{\Q^{ra}} \mathrm{HS}^{ra}({K}_L^\perp \1_{x})(x)$ plus a remainder to see
				\begin{equation} \label{e.localize_1} 
					\sum_{x \in \car_{L - \ell}} a(x)\inner{a, \Esp_{\Q^{ra}} \mathrm{HS}^{ra}{K}_L^\perp \1_{x}}_{\car_L} = \sum_{x \in \car_{L - \ell}} |a(x)|^2 \(\Esp_{\Q^{ra}} \mathrm{HS}^{ra}{K}_L^\perp \1_{x}\)(x) + \mathrm{Error}_1
				\end{equation}
				where
				\begin{align*}
					\mathrm{Error}_1 &\leq C \sum_{x \in \car_{L- \ell}} |a(x)| \| a \|_{L^2(\car_L)} \| \Esp_{\Q^{ra}} \mathrm{HS}^{ra} K_L^\perp \1_x \|_{L^2(\car_L \setminus \{x\})} \\ &\leq C \sum_{x \in \car_{L - \ell}} \( |a(x)| \| a \|_{L^2(\car_L)} \left [\| {K}_L^\perp \1_x \|_{L^2(\car_L \setminus \{ x \})} + \sqrt{B_L(K_L^\perp \1_x)} \right ] \).
				\end{align*}
				We used \eref{ptmass} (with $u_a$ in \eref{ptmass} replaced by $u_{K_L^\perp \1_x} = \mathrm{HS}^{ra} K_L^\perp \1_x$) in passing to the last line above. Since ${K}_L^\perp \1_x(y) = -{K}_L \1_x(y)$ for $y \ne x$, \eref{berg_infvol} proves the bound
				\begin{align*}
					\lefteqn{\bigg | \sum_{x \in \car_{L - \ell}} |a(x)| \| a \|_{L^2(\car_L)} \| {K}_L^\perp \1_x \|_{L^2(\car_L \setminus \{ x \})} \bigg |} \quad & \\ &\leq C \| a \|_{L^2(\car_L)} \| a \|_{L^1(\car_L)} \ell^{-d/2} \leq C L^{-d/2 + 2 + \epsilon d/2}
				\end{align*}
				which controls part of $\mathrm{Error}_1$. We now control $B_L(K_L^\perp \1_x)$. By definition, we have 
				$$
				\inner{K_L^\perp \1_x, P^z}_{\cl_1 \car_L} = 0
				$$
				for all $z \in \pa \car_L$. Thus
				$$
				\inner{K_L^\perp \1_x, P^z}_{\car_L} = \inner{K_L^\perp \1_x, P^z}_{\cl_1 \car_L} - K_L^\perp \1_x(z) = -K_L^\perp \1_x(z) = K_L \1_x(z).
				$$
				Thus we can again apply \eref{berg_infvol} to get
				\begin{equation}\label{e.BLKperp1}
					B_L(K_L^\perp \1_x) := \sum_{z \in \pa \car_L} \bigg |\inner{K_L^\perp \1_x, P^z}_{\car_L}  \bigg|^2 \leq \| K_L \1_x \|_{L^2(\cl_1 \car_L)}^2 \leq C \rho_x^{-d}.
				\end{equation}
				We use this to get a final bound
				$$
				\mathrm{Error}_1 \leq C L^{-d/2 + 2 + \epsilon d/2}.
				$$
				
				It remains to further simplify the RHS of \eref{localize_1}. By \eref{meanthermo}, we have
				\begin{equation} \label{e.localize_2}
					\sum_{x \in \car_{L - \ell}} |a(x)|^2 \(\Esp_{\Q^{ra}} \mathrm{HS}^{ra}{K}_L^\perp \1_{x}\)(x) = \sum_{x \in \car_{L - \ell}} |a(x)|^2 \Esp_{\nu^{ra(x)}} U^{ra(x)} K_L^\perp \1_x(x) + \mathrm{Error}_2
				\end{equation}
				where
				$$
				|\mathrm{Error}_2| \leq  C\sum_{x \in \car_{L- \ell}} |a(x)|^2\(|{K}_L^\perp \1_x(x)| \ell^{-d/2} + \| {K}_L^\perp \1_x \|_{L^2(\car_L \setminus \{x\})} + \sqrt{B_L(K^\perp_L \1_x)}\).
				$$
				The term $\mathrm{Error}_2$ is bounded by $C \ell^{-d/2}$ by the same method as the bound for $\mathrm{Error}_1$, i.e.\ by \eref{berg_infvol} and \eref{BLKperp1}. We can replace $K_L^\perp \1_x(x)$ by $1$ in the RHS of \eref{localize_2} by using the estimate
				$$
				|K_L^\perp \1_x(x) - 1| \leq \| K_L \1_x \|_{L^2(\cl_1 \car_L)} \leq C \rho_x^{-d/2},
				$$
				and the resulting error is bounded by $C \ell^{-d/2}$.
				
				Combining our bounds on $T_1$ and $T_2$ with \eref{localize_1} and \eref{localize_2}, we have almost proved \eref{thermoform}, except our representation involves summing $|a(x)|^2 \Esp_{\nu^{ra(x)}} U^{ra(x)}$ over $x \in \car_{L - \ell}$ instead of $\car_L$. This is fixed by noting
				$$
				\bigg | \sum_{x \in \Lambda} |a(x)|^2 \Esp_{\nu^{ra(x)}} U^{ra(x)} \bigg | \leq C \| a \|_{L^2(\Lambda)}^2 \leq CL^{-d + 4 - \epsilon},
				$$
				and \eref{thermoform} is proved.
				
				Next, we explain how to conclude \tref{inf_limit_thm}. First, we translate \eref{inf_limit} to the level of $\Q^0$ and apply Helffer-Sj\"ostrand representation in \eref{cumulant_var} and \eref{var_form}:
				$$
				\log \Esp_{\MM_L} \exp\( \sum_{x \in \car_K} a'(x) \phi(x) \) = \log \Esp_{\Q^0} e^{ \inner{a, \eta}_{\car_L}} = \int_0^1 (1-r) \inner{a, \Esp_{\Q^{ra}} \mathrm{HS}^{ra} a}_{\car_L} dr.
				$$
				We apply \eref{thermoform} to get
				$$
				\log \Esp_{\MM_L} \exp\( \sum_{x \in \car_K} a'(x) \phi(x) \) = \sum_{x \in \car_L} |a(x)|^2 \int_0^1 (1-r)\Esp_{\nu^{ra(x)}} U^{ra(x)} dr + \mathrm{Error}
				$$
				for $|\mathrm{Error}| \leq CL^{-d/2 + 2 + \epsilon d/2}$.
				The Helffer-Sj\"ostrand representation associated to the measure $\nu^{ra(x)}$, analogous to \eref{cumulant_var} and \eref{var_form}, is
				$$
				\Var_{\nu^{ra(x)}} [\xi] = \Esp_{\nu^{ra(x)}} U^{ra(x)}, \quad \log \Esp_{\nu^0} \exp(a(x) \xi) = \int_0^1 (1-r) \Var_{\nu^{ra(x)}} [a(x) \xi] dr,
				$$
				and so
				\begin{equation} \label{e.inf_vol_rep1}
					\log \Esp_{\MM_L} \exp\( \sum_{x \in \car_K} a'(x) \phi(x) \) = \sum_{x \in \car_L} \log \Esp_{\nu^0} \exp( (\Gamma_L \ast a')(x) \xi) + \text{Error}.
				\end{equation}
				All that remains is to pass to $L \to \infty$ in the RHS. Let $\psi(\lambda) = \log \Esp_{\nu^0} e^{\lambda \xi}$ be the c.g.f.\ of $\nu^0$, which has locally bounded derivatives. We have
				$$
				|\psi(\Gamma_L \ast a'(x)) - \psi(\Gamma \ast a'(x))| \leq C |\Gamma_L \ast a'(x) - \Gamma \ast a'(x)|
				$$
				where $C$ depends only on $a'$. Using the estimates in \pref{greens_theo} and \eref{gam_rep}, we can bound the RHS by
				\begin{align*}
					|\Gamma_L \ast a'(x) - \Gamma \ast a'(x)| &\leq \sum_{z \in \pa \car_L} \sum_{y \in \car_K} P^z_L(x)\Gamma(z,y) |a'(y)| \\ &\leq C \| a' \|_{L^1(\car_K)} \sup_{z \in \pa \car_L, y \in \car_K} \frac{1}{|z - y|^{d-2}} \leq C L^{-d+2}.
				\end{align*}
				Observe that $\Gamma_L \ast a'(x)$ decays like $|x|^{-d + 2}$ as $|x| \to \infty$, uniformly in $L$, and $\Esp_{\nu^0} \xi = 0$. It follow that the c.g.f.\ of $\nu^0$ at $\Gamma_L \ast a'(x)$ decays like $|x|^{-2d+4}$, and so, for any fixed $\beta > 0$, we have
				$$
				\sum_{x \in \Z^d \setminus \car_{L^\beta}} |\psi(\Gamma_L \ast a'(x))| \leq C \sum_{x \in \Z^d \setminus \car_{L^\beta}} \frac{1}{|x|^{2d-4}} = O(L^{-\beta(d - 4)})
				$$
				as $L \to \infty$. The same estimate holds for $\Gamma$ in place of $\Gamma_L$. Combining the above, we see
				\begin{align*}
					\sum_{x \in \car_L} \psi (\Gamma_L \ast a'(x)) &= \sum_{x \in \car_{L^\beta}} \psi(\Gamma_L \ast a'(x)) + O(L^{-\beta(d-4)}) \\ &=  \sum_{x \in \car_{L^\beta}} \psi(\Gamma \ast a'(x)) + O(L^{-d+2} L^{\beta d}) + O(L^{-\beta(d-4)}) \\
					&= \sum_{x \in \Z^d} \psi(\Gamma \ast a'(x)) + O(L^{-d+2} L^{\beta d}) + O(L^{-\beta(d-4)}).
				\end{align*}
				Choosing $\beta > 0$ sufficiently small finishes the proof of \tref{inf_limit_thm}.
			\end{proof}
		\end{subsection}

		\begin{subsection}{Limits of the rescaled model} Next, we prove the scaling limit in $d \geq 2$. Let $f : [-1,1]^d \to \R$ and let $u$ be the solution to continuum problem
			\begin{equation}\label{e.cont_bil_eq}
				\begin{cases}
					{\Delta_{\R^d}^2} u(x) = f(x), \quad &x \in (-1,1)^d,\\
					u(x) = \pa_n u(x) = 0, \quad &x \in \pa (-1,1)^d.
				\end{cases}
			\end{equation}
			Depending on the dimension $d$, we will make one of the two following assumptions:
			\begin{equation} \label{e.A1}
				u \in C^5([-1,1]^d) \text { with } \| u \|_{C^5([-1,1]^d)} \leq C < \infty
			\end{equation}
			or
			\begin{equation} \label{e.A2}
				f \in C^0([-1,1]^d) \text{ with } \| f \|_{L^\infty([-1,1]^d)} \leq C < \infty.
			\end{equation}
			Condition \eref{A2} is clearly weaker than \eref{A1}, and it suffices in the case of square/cubic domains in $d=2,3$ using estimates in \cite{MS19}. We conjecture that it is also sufficient in all other cases. We set $a(x) = -L^{-d/2 - 2} \sum_{y \in \car_L} f(y/L) \Gamma_L(x,y)$.
			
			We now cite some basic properties of $a$ and ${K}_L$ proved in \sref{append}. The field ${K}_L a$ is harmonic in $\car_L$. Let $\Lambda' = \cl_1 \car_L \setminus \car_{L - \ell'}$ for some integer $\ell' \in [1, L]$, i.e.\ $\Lambda'$ is a boundary layer of $\cl_1 \car_L$. If we assume \eref{A1} holds or if we assume $d=2,3$ and \eref{A2} holds, then by either \eref{olKbdry} or \tref{DeltaGthm}, we have
			\begin{equation} \label{e.KLa_bdry_2}
				\| {K}_L a \|_{L^2(\Lambda')}^2 + \| (\text{Id} - {K}_L)a \|_{L^2(\Lambda')}^2 \leq \frac{C \cdot (\ell'+1)}{L}(1 + \1_{d=2} (\log L)^2).
			\end{equation}
			The same inequality holds for $\Lambda' = \cl_1 \car_L$ if we set $\ell' = L$. It is important that $\Lambda'$ contains $\pa \car_L$, unlike the layer $\Lambda$ considered in the previous subsection, and the constant $C$ depends only on the dimension and the size of the constants in \eref{A1} or \eref{A2}. We also have
			\begin{equation} \label{e.ainfty1}
				\| a \|_{L^\infty(\cl_1 \car_L)} \leq CL^{-d/2} \| f \|_{L^\infty([-1,1]^d)} (1 + \1_{d=2} \log L) \leq CL^{-d/2}(1 + \1_{d=2} \log L),
			\end{equation}
			which follows easily from a standard estimate on $\Gamma_L$ stated in \eref{gam_sup}.
			
			\begin{theo} \label{t.scal_lim}
				Let $a$ be as above and assume either \eref{A2} holds and $d = 2,3$ or \eref{A1} holds and $d \geq 2$. Then for any $r \in [0,1]$ we have
				\begin{equation} \label{e.scal_lim_int_est}
					\inner{a, \Esp_{\Q^{ra}} \mathrm{HS}^{ra} a}_{\car_L} = \Esp_{\nu^0}U^0 \| (\mathrm{Id} - {K}_L) a \|_{L^2(\cl_1 \car_L)}^2 + \mathrm{Error},
				\end{equation}
				where $\mathrm{Error}$ satisfies
				$$
				|\mathrm{Error}| \leq CL^{-\frac{d-1}{6d-2}}(1 + \1_{d=2}(\log L)^3).
				$$
				As a consequence, \tref{scal_limit_int} holds.
			\end{theo}
			\begin{proof}
				We first prove \eref{scal_lim_int_est} and then explain how it implies \tref{scal_limit_int}.
				
				We let $\epsilon > 0$ be a parameter to be fixed later. Let $\Lambda = \car_L \setminus \car_{L - \ell}$ be a boundary layer of width $\ell$ between $L^{1-\epsilon}$ and $2L^{1-\epsilon}$. Let $\Lambda' = \Lambda \cup \pa \car_L$.
				
				Our strategy for approximating $\Esp_{\Q^{ra}} \mathrm{HS}^{ra} a$ is to decompose $a$ as ${K}_L a + {K}_L^\perp a$ and compute the linear operator $\Esp_{\Q^{ra}} \mathrm{HS}^{ra}$ on each piece.
				
				The first piece ${K}_L a$ is a harmonic function which does not concentrate on $\pa \car_L$ by \eref{KLa_bdry_2}. The estimate \eref{Estiff1} from \pref{Estiffprop} therefore applies, and we have
				\begin{align*}
					\| \Esp_{\Q^{ra}} \mathrm{HS}^{ra} {K}_L a \|_{L^2(\car_L)}^2 &\leq C \mathcal{E}_0(\mathrm{HS}^{ra}{K}_L a; \Q^{ra}) = C\mathcal{E}_0(u_{{K}_L a}; \Q^{ra})\\  &\leq C \| {K}_L a \|_{L^2(\car_L)} \| {K}_L a \|_{L^2(\pa \car_L)} \leq CL^{-1/2}(1 + \1_{d=2} (\log L)^2).
				\end{align*}
				In the last line, we used \eref{KLa_bdry_2} on the layer $\Lambda' = \pa \car_L$ of width $1$ to estimate the boundary term.
				
				Next, we consider the term coming from ${K}_L^\perp a$. Applying the result in equation \eref{mean_scaling} gives
				\begin{align} \label{e.scal_lim_perp}
					\lefteqn{ \| \Esp_{\Q^{ra}} \mathrm{HS}^{ra} {K}_L^\perp a - \Esp_{\nu^0} U^0  \cdot {K}_L^\perp a \|_{L^2(\car_L)}} \quad & \\\notag &\leq CL^{-d/2 + 1/2 + c_d \epsilon}(1 + (\log L)^2 \1_{d=2}) + C \| {K}_L^\perp a \|_{L^2(\Lambda)} \\ \notag &\quad + C\sqrt{B_L(K_L^\perp a)} + C \( \sum_{x \in \car_L} |a(x)|^2 |{K}_L^\perp a(x)|^2 \)^{1/2},
				\end{align}
				where $c_d = \frac32 d - 1$ and $B_L$ is defined in \eref{BL_def}. We rewrite $B_L(K_L^\perp a)$ like in \eref{BLKperp1} to get
				$$
				B_L(K_L^\perp a) = \sum_{z \in \pa \car_L} \bigg | \inner{K_L^\perp a, P^z}_{\car_L} \bigg |^2 = \| K_L^\perp a \|_{L^2(\pa \car_L)}^2.
				$$
				Considering the estimates \eref{KLa_bdry_2} and \eref{ainfty1} on $a$, it is straightforward to bound the RHS of \eref{scal_lim_perp} by
				$$
				C(L^{-d/2 + 1/2 + c_d \epsilon} + L^{-\frac12 \epsilon} + L^{-d/2})(1 + (\log L)^2 \1_{d=2}).
				$$
				We put together the preceding estimates along with Cauchy-Schwarz to conclude
				$$
				\inner{a, \Esp_{\Q^{ra}} \mathrm{HS}^{ra}a}_{\car_L} = \Esp_{\nu^0}U^0 \inner{a,  {K}_L^\perp a}_{\car_L } + \text{Error}
				$$
				where $\text{Error}$ satisfies
				$$
				|\text{Error}| \leq C \(L^{-d/2 + 1/2 + (3d/2 - 1) \epsilon} + L^{- \epsilon/2} + L^{-d/2} + L^{-1/4}\)(1 + \1_{d=2} (\log L)^3).
				$$
				We optimize by choosing $\epsilon = \frac{d-1}{3d-1}$, which gives
				$$
				|\text{Error}| \leq CL^{-\frac{d-1}{6d-2}}(1 + \1_{d=2}(\log L)^3).
				$$
				Finally, since ${K}_L$ is an $L^2(\cl_1 \car_L)$ projection and $a(z) = 0$ for $z \in \pa \car_L$, we have
				$$
				\inner{a, {K}_L^\perp a}_{\car_L } = \inner{a, {K}_L^\perp a}_{\cl_1 \car_L } = \|{K}_L^\perp a\|_{L^2(\cl_1 \car_L)}^2,
				$$
				which concludes the proof of \eref{scal_lim_int_est}.
				
				To conclude \tref{scal_limit_int}, notice that $\Esp_{\nu^0} U^0$ is equal to the variance of a variable $\xi$ with law $\nu^0$, and the approximation \eref{scal_lim_int_est} can also be done for a Gaussian membrane model. Then the representations \eref{cumulant_var} and \eref{var_form} relate \eref{scal_lim_int_est} to the cumulant generating function, and allow us to conclude the theorem by comparing the approximations for the Gaussian and non-Gaussian model.
			\end{proof}
			
			Our theory makes no serious distinction between the first estimate of \tref{scal_limit_int} and \tref{scal_limit_pt}, and we now prove the latter.
			\begin{proof}[Proof of \tref{scal_limit_pt}]
				The proof is the same as that of \tref{scal_lim}, except with a different form for $a$, which is $a(x) = L^{d/2 -2}\sum_{i=1}^k \Gamma_L(x,L  y_i) c_i$ for some numbers $c_i$ and points $y_i \in [-1,1]^d$. The values $c_i$ and $y_i$ change with $L$ because of the continuous interpolation of the interface $\ol{\phi}$, but we have effectively the same estimates on $a$ as \eref{KLa_bdry_2}. In particular, we have
				$$
				\| a \|_{L^\infty(\car_L)} \leq CL^{d/2-2}(1+ \1_{d=2} \log L), \quad \| a \|_{L^2(\car_L)} \leq C(1 + \1_{d=2} \log L)
				$$
				uniformly over $L$. And for a boundary layer $\Lambda'$ of width $\ell$, using \tref{DeltaGthm}, we have
				\begin{equation} \begin{aligned}
						\| {K}_L^\perp a(x) \|_{L^2(\Lambda')}^2 &\leq C(1 + \1_{d=2}(\log L)^2) L^{d-4} \sup_i \sum_{x \in \Lambda'} \frac{1}{(1 + |x-Ly_i|)^{2d-4}} \\
						&\leq  C(1 + \1_{d=2}(\log L)^2) L^{d-4} \ell \sum_{x \in \Z^{d-1}, \| x \| \leq L} \frac{1}{(1 + |x|)^{2d-4}} \\
						&\leq C(1 + \1_{d=2}(\log L)^2) \frac{\ell}{L}.
				\end{aligned} \end{equation}
				On the middle line above, we dominated the sum over the boundary layer by $\ell$ times a sum over a $d-1$ dimensional plane of length $L$. By the same idea, we have 
				$$
				\| a \|_{L^2(\Lambda')}^2 \leq CL^{-1} \ell (1 + \1_{d=2} (\log L)^2).
				$$
				
				The bounds above are exactly the same in $d=2$ as those used in the proof of \tref{scal_lim}. For $d=3$, we have a weaker $L^\infty$ control on $a$ than before, but it is inconsequential for the final error estimate. We omit the rest of the proof the main estimate.
				
				In order to prove the existence of a Gaussian distributional limit in the space of continuous functions on $[-1,1]^d$, we must address tightness of the distribution of $\ol{\phi}$ on $C([-1,1]^d)$. 
				
				Let $x,y \in [-1,1]^d$. We wish to bound from above the typical size of $|\ol{\phi}(x) - \ol{\phi}(y)|$ in terms of $|x-y|$. The result we desire is proved in \cite{CDH19} for the Gaussian membrane model, so we seek to use Gaussian domination coming from $\inf V'' > 0$. To do this, we need a version of the Brascamp-Lieb inequality from \cite{BL76}, Theorem 5.1, which is a generalization of the Poincar\'e inequality \pref{poincarebe}. The Hessian of the general membrane model energy is larger than the Hessian of the Gaussian membrane model associated to the single-spin potential $\xi \mapsto C^{-1} \xi^2$, where $C$ depends on $\inf V''$. We have
				$$
				\Esp_{\MM_L} \left [ |\ol{\phi}(x) - \ol{\phi}(y) |^{2\beta}\right] \leq C \Esp_{\MM^G_L}\left[|\ol{\phi}(x) - \ol{\phi}(y)] |^{2\beta} \right]
				$$
				where $\MM^G_L$ is a Gaussian membrane model and $2\beta \geq 1$. It follows from Lemma 2.6 in \cite{CDH19} that
				$$
				\Var_{\MM_L^G} [\ol{\phi}(x) - \ol{\phi}(y)] \leq C |x-y|^{1 + \lambda}
				$$
				for any fixed $\lambda \in (0,1)$ in $d=2$ or $\lambda = 0$ in $d=3$. Since the model is Gaussian, it follows
				$$
				\Esp_{\MM^G_L}\left[|\ol{\phi}(x) - \ol{\phi}(y) |^{2\beta} \right] \leq C |x-y|^{\beta(1+\lambda)},
				$$
				and so we conclude the same inequality for $\MM_L$, except with a potentially larger constant $C$.
				
				Tightness and the H\"older continuity properties of the limit follow from the Kolmogorov-Chentsov criterion (see e.g.\ \cite{K02}, Corollary 16.9). Since \cite{CDH19} proves the limit of $\ol{\phi}$ under $\MM^G_L$, and we know these limits must agree up to a scaling factor, the proof is concluded.
			\end{proof}
		\end{subsection}
	\end{section}

	\begin{section}{Bergman projection and Green's function estimates} \label{s.append}
		In this section, we collect already known estimates for the Laplacian Green's functions $-\Gamma_L : \car_L \times \car_L \to \R$. We also prove a non-concentration result for the discrete harmonic Bergman projection. Many of the estimates are not sharp, but we include them to be self-contained and to show that, for the most part, only weak bounds are needed for our results. We also could not find a reference for the discrete harmonic Bergman projection results.
		
		\begin{subsection}{Laplacian Green's function Estimates}
			The first result gives some non-optimal bounds on the Poisson kernel and Green's function of a simple random walk. We have chosen to give weak bounds because the proof is easy to adapt to domains other than $\car_L$.
			\begin{prop} \label{p.greens_theo}
				Let $d \geq 2$, and let $-\Gamma_L$ be the Dirichlet Green's function for the discrete Laplacian on $\car_L$. Then there is a constant $C$, depending only on $d$, such that
				\begin{equation} \label{e.gam_sup}
					|\Gamma_L(x,y)| \leq C (1 + |x-y|)^{2-d}(1 + \1_{d=2} \log L)
				\end{equation}
				for all $x,y \in \car_L$. We also have (recall that $\rho_y = \dist(y, \pa \car_L)$)
				\begin{equation} \label{e.greens_deriv}
					|\nab_y \Gamma_L(x,y)| \leq C (1 + |x-y|)^{1-d} + C \rho_y^{1-d}
				\end{equation}
				where $\nab_y$ denotes the discrete derivative in $y$. We can use this to bound the Poisson kernel $P^z_L \in \R^{\car_L}$:
				\begin{equation} \label{e.poisson}
					|P^z_L(y)| \leq C \rho_y^{-d+1}.
				\end{equation}
				
			\end{prop}
			
			\begin{proof}
				We consider dimensions $d \geq 3$ first. Let $\Gamma$ be the Green's function for $-\Delta$ on $\Z^d$. We have the representation
				\begin{equation} \label{e.gam_rep}
					\Gamma_L(x,y) = \Gamma(x,y) - \sum_{z \in \pa \car_L} P^z_L(x) \Gamma(z,y),
				\end{equation}
				which can be proved by noticing the RHS solves the same discrete Poisson equation as $\Gamma_L$. Estimates \eref{gam_sup} and \eref{greens_deriv} are true with $\Gamma$ in place of $\Gamma_L$, as proved in \cite{L91}. Estimate \eref{gam_sup} follows then directly from the positivity of $\Gamma_L$ and $\Gamma$. Taking the derivative of the representation and applying estimates for $\nab_y \Gamma$, we can write
				\begin{align}
					|\nab_y \Gamma_L(x,y)| &\leq |\nab_y \Gamma(x,y)| + \sum_{z \in \pa \car_L} P^z_L(x) |\nab_y \Gamma(z,y)| \\ \notag
					&\leq C (1 + |x-y|)^{1-d} + C\rho_y^{-d+1} \sum_{z \in \pa \car_L} P^z_L(x) \\ \notag &= C (1 + |x-y|)^{1-d} + C\rho_y^{-d+1}.
				\end{align}
				In the last inequality, we used the fact that $\sum_z P^z_L(x) = 1$ is the probability that a simple random walk starting at $x$ exits $\car_L$ in finite time.
				
				For $d=2$, the infinite volume Green's function does not exist, but the potential kernel $a(y-x) = \lim_{L \to \infty} (\Gamma_L(x,y) - \Gamma_L(x,x))$ exists and depends only on $y-x$. The representation formula \eref{gam_rep} holds in the sense
				$$
				\Gamma_L(x,y) = \bigg(\sum_{z \in \pa \car_L} P^z_L(x) a(z-y) \bigg) - a(y-x).
				$$
				The estimates \eref{gam_sup} and \eref{greens_deriv} follow in the same manner as before, now substituting the estimates $|a(x)| \leq C \log |x|$ and $|\nab_x a(x)| \leq C|x|^{-1}$.
				
				Note that $P^z_L(y) = \Gamma_L(y, \tilde{z})$, where $\tilde{z} \in \car_L$ is the unique point adjacent to $z \in \pa \car_L$. Indeed, the function $y \mapsto \Gamma_L(y, \tilde{z}) + \1_z$ is harmonic in $\car_L$ with boundary condition $\1_z$, and so it must coincide with $P^z_L$. The estimate \eref{poisson} follows from \eref{greens_deriv} by noting $\Gamma_L(y,z) = 0$.
			\end{proof}	
		\end{subsection}
		
		\begin{subsection}{The harmonic Bergman projection}
			In this subsection, we give estimates for the Bergman projection $K_L$ applied to $a$ of the form that we encounter in the proofs of our main theorems.
			
			Let $\mathcal{H}(\car_L)$ denote the subspace of functions on $\car_L \cup \pa \car_L$ which are discrete harmonic in $\car_L$. The Bergman projection $K_L : \R^{\car_L \cup \pa \car_L} \to \mathcal{H}(\car_L)$ is defined by
			\begin{equation} \label{e.KL_def}
				K_L a = \text{argmin}_{\tilde{a} \in \mathcal{H}(\car_L)} \| \tilde{a} - a \|_{L^2(\car_L \cup \pa \car_L)}^2.
			\end{equation}
			
			We can relate $K_L$ to the bi-Laplacian Green's function $G_L$, which is defined as the solution to
			$$\begin{cases}
				\Delta^2_x G_L(x,y) = \1_{x=y}, \quad &x \in \car_L,\\
				G_L(x,y) = 0, \quad &x \in \Z^d \setminus \car_L.
			\end{cases} $$
			By a slight abuse of notation, we can consider $G_L$ as a convolutional operator $\R^{\car_L} \to \R^{\car_L \cup \pa^2 \car_L}$. Recall that $\pa^2 \car_L$ consists of all points in $\Z^d \setminus \car_L$ within $\ell^1$ distance $2$ of $\car_L$. We can also consider the Laplacian $\Delta$ as an operator $\R^{\car_L \cup \pa \car_L} \to \R^{\car_L}$ or $\R^{\car_L \cup \pa^2 \car_L} \to \R^{\car_L \cup \pa \car_L}$. Under these identifications, we now verify that $K^\perp_L = \Delta G_L \Delta$, where $K^\perp_L = \text{Id} - K_L$.
			
			Indeed, for any $v$ in the image of $K^\perp_L$, we have that $v$ is $L^2(\car_L \cup \pa \car_L)$-orthogonal to the kernel of $\Delta : \R^{\car_L \cup \pa \car_L} \to \R^{\car_L}$, and so is in the image of the adjoint $\Delta : \R^{\car_L} \to \R^{\car_L \cup \pa \car_L}$ given by extension by $0$ outside $\car_L$. That is, $v = \Delta u$ for some $u \in \R^{\Z^d}$ with $\supp (u) \subset \car_L$. And so $\Delta G_L \Delta v = \Delta G_L \Delta^2 u = \Delta u = v$. We also clearly have $(\Delta G_L \Delta)(\mathcal{H}(\car_L)) = \{ 0 \}$, and so the claim $K^\perp_L = \Delta G_L \Delta$ is established.
			
			For our main theorems, we need to estimate $K_L a$ for $a = \Gamma_{L}\ast f_L$ where $f_L$ is a "nice" element of $\R^{\cl_1 \car_L}$. Specifically, we need to prove that the $L^2$ mass of $K_L a$ does not concentrate in a boundary layer of $\car_L$.
			
			One convenient condition on $f_L$ is that it is a rescaled version of a continuum function $f$ for which the continuum bi-Laplacian problem has regular solution. This is assumption \eref{A1} in \sref{concl_sec}.
			
			For $f : [-1,1]^d \to \car_L$, define $f_L(x) = L^{-d/2 - 2} f(x/L)$ as a mapping $\car_L \to \R$. Let $u : [-1,1]^d \to \car_L$ solve
			\begin{equation}
				\begin{cases}
					{\Delta}_{\R^d}^2 u(x) = f(x), \quad &x \in (-1,1)^d,\\
					u(x) = \pa_n u(x) = 0, \quad &x \in \pa (-1,1)^d,
				\end{cases}
			\end{equation}
			where $\pa_n$ denotes the normal derivative and ${\Delta}_{\R^d}$ denotes the continuum Laplacian. Define $u_{L} :  \car_L \cup \pa^2 \car_L \to \R$ by $u_L(x) = L^{-d/2 + 2} u(x/(L+2))$. Finally, let $v : \cl_2 \car_L \to \R$ be defined by
			\begin{equation}
				\begin{cases}
					\Delta^2 v(x) = f_L(x), \quad &x \in \car_L,\\
					v(x) = 0, \quad &x \in \pa^2 \car_L,
				\end{cases}
			\end{equation}
			where $\Delta^2$ is the discrete bi-Laplacian operator.
			
			The following proposition shows that $\Delta v = \Delta G_L f_L$ and $\Delta u_L$ are close in $L^2$ assuming $u$ is regular. Note that $K_L = \text{Id} - \Delta G_L \Delta$, and, using the below proposition, we may estimate $K_L \Gamma_{L} \ast f_L = \Gamma_L\ast f_L - \Delta v$ by comparing to a continuum object.
			\begin{prop}
				Let $u, u_L, v \in \R^{\cl_2 \car_L}$ be defined as above, and assume $u$ is five times differentiable in $(-1,1)^d$. Then we have
				\begin{equation} \label{e.bil_cont_approx}
					\| \Delta u_L - \Delta v \|_{L^2(\cl_1 \car_L)} \leq CL^{-1/2} \| u \|_{C^5([-1,1]^d)}
				\end{equation}
				for a constant $C$ dependent only on $d$.
			\end{prop}
			\begin{proof}
				Let $w = u_L - v$. The proof will follow two steps: (1) we can bound the quantity $\| \Delta w \|_{L^2(\cl_1 \car_L)}$ by $\| \Delta^2 w \|_{L^2(\car_L)}$ and the boundary values of $w$ on $\pa^2 \car_L$, and (2) $w$ is almost a biharmonic function with Dirichlet boundary conditions.
				
				By discrete integration by parts, we have
				$$
				\| \Delta w \|_{L^2(\cl_1 \car_L)}^2 = \inner{w, \Delta^2 w}_{\car_L} + \sum_{x \in \pa^2 \car_L, y \in \cl_1 \car_L} \Delta(x,y) w(x) \Delta w(y)
				$$
				where $\Delta(x,y)$ is $1$ if $x \sim y$ and $-2d$ if $x = y$ and $0$ otherwise. Thus
				\begin{equation} \label{e.wparts1}\begin{aligned}
						\| \Delta w \|_{L^2(\cl_1 \car_L)}^2 &\leq \| w \|_{L^2(\car_L)} \| \Delta^2 w \|_{L^2(\car_L)} + C \| w \|_{L^\infty(\pa^2 \car_L)} \sum_{y \in \cl_1 \car_L\setminus \car_{L-1}} |\Delta w(y)| \\
						&\leq \| w \|_{L^2(\car_L)} \| \Delta^2 w \|_{L^2(\car_L)} + C L^{(d-1)/2}  \| w \|_{L^\infty(\pa^2 \car_L)} \| \Delta w \|_{L^2(\cl_1 \car_L)}.
				\end{aligned}\end{equation}
				We extend $w$ by $0$ to $\car_{L+10}$ (say), let $\nab^2 w$ denote the matrix of second discrete derivatives of $w$, and apply the discrete Poincar\'e inequality twice to get
				\begin{equation} \label{e.wpoinc1} \begin{aligned}
						\| w \|_{L^2(\car_L)} \leq CL^2 \| \nab^2 w \|_{L^2 (\car_{L+10})} &= CL^2 \| \Delta w \|_{L^2(\car_{L+10})} \\ &\leq C L^2 \| \Delta w \|_{L^2(\cl_1 \car_L)} + CL^{2 + (d-1)/2} \| w \|_{L^\infty(\pa ^2 \car_L)}.
				\end{aligned} \end{equation}
				The last inequality follows from the fact that $|\Delta w(x)| \leq 4d \max_{y, \dist(y,x) \leq 1} |w(y)|$.
				Applying Young's inequality to \eref{wparts1} gives
				\begin{align*}
					\lefteqn{\| \Delta w \|_{L^2(\cl_1 \car_L)}^2} \quad & \\ &\leq \frac{\epsilon^2}{L^4} \| w \|_{L^2(\car_L)}^2 + C_\epsilon L^4 \| \Delta^2 w \|_{L^2(\car_L)} + \epsilon^2 \| \Delta w\|_{L^2(\cl_1(\car_L))}^2 + C_\epsilon L^{d-1} \| w \|_{L^\infty(\pa^2 \car_L)}^2.
				\end{align*}
				for any $\epsilon > 0$. Choosing $\epsilon^2 < 1/C$ and using \eref{wpoinc1} shows
				\begin{equation} \label{e.deltaw_est1}
					\| \Delta w \|_{L^2(\cl_1 \car_L)}^2 \leq CL^4 \| \Delta^2 w \|_{L^2(\car_L)}^2 + CL^{d-1} \| w \|_{L^\infty(\pa^2 \car_L)}^2,
				\end{equation}
				which completes step (1) of the proof.
				
				We now bound $\Delta^2 w$ in $\car_L$. Let $L' = L + 2$ to lighten notation. Applying repeatedly the fundamental theorem of calculus, we compute
				\begin{align}
					\Delta u_L(x) &= L^{-d/2 + 2} (L')^{-2} {\Delta}_{\R^d} u(x/L') \\\notag & \quad + L^{-d/2 + 2} (L')^{-2}\sum_{i=1}^d \sum_{\sigma = \pm 1} \int_0^1 (1-s) s \( \nab^2_{i,i} u\(\frac{x + \sigma s e_i}{L'}\) - \nab^2_{i,i} u\(\frac{x}{L'}\) \) ds.
				\end{align}
				where $\nabla^2_{i,i}$ denotes the continuum second derivative in the standard basis direction $e_i$. The term on the last line can be bounded by $CL^{-d/2 - 1} \| \nab_{\R^d}^3 u \|_{L^\infty}$. Moreover, a similar argument applied to the integrand shows the last term has second discrete derivative bounded by $CL^{-d/2 - 3} \| \nab_{\R^d}^5 u \|_{L^\infty}$. So we can iterate the above identity to see
				$$
				\Delta^2 u_L(x) = L^{-d/2 - 2} (\Delta_{\R^d}^2 u)(x/L) + \text{Error}_x = f_L(x) + \text{Error}_x
				$$
				where
				$$
				|\text{Error}_x|\leq C  L^{-d/2 - 3} \| \nab_{\R^d}^5 u \|_{L^\infty} + C L^{-d/2 - 3} \| \Delta_{\R^d}^2 u \|_{L^\infty} + C L^{-d/2 -3} \| \nab_{\R^d} f \|_{L^\infty}
				$$
				for a constant $C$ independent of $L$ and $u$. Above, we have also changed $L'$ into $L$ and generated the corresponding error terms. Using $\Delta^2 v = f_L$, we conclude that
				$$
				\| \Delta^2 w \|_{L^2(\car_L)} \leq \| \text{Error}_{\cdot} \|_{L^2(\car_L)} \leq C L^{-3} \( \| \nab_{\R^d}^5 u \|_{L^\infty} + \| \Delta_{\R^d}^2 u \|_{L^\infty} +\| \nab_{\R^d} f \|_{L^\infty} \).
				$$
				
				Next, since $w(x) = L^{-d/2 + 2}u(x/(L+2))$ for $x \in \pa^2 \car_L$ and by the boundary conditions of $u$, we have
				$$
				\| w \|_{L^\infty(\pa^2 \car_L)} \leq C L^{-d/2} \| \nab_{\R^d}^2 u \|_{L^\infty}.
				$$
				This completes step (2) of the proof. We conclude from the estimates on $w$ and $\Delta^2 w$ and \eref{deltaw_est1} that 
				$$
				\| \Delta u_L - \Delta v \|_{L^2(\cl_1 \car_L)} \leq C L^{-1/2}\| u \|_{C^5([-1,1]^d)},
				$$
				and the proof is complete.
			\end{proof}
			
			We now bound the Bergman projection $K_L$ on a boundary layer.
			\begin{prop}
				Let $a = \Gamma_L \ast f_{L} \in \R^{\cl_1 \car_L}$ and $f_L = L^{-d/2 - 2} f(x/L)$ for a function $f : (-1,1)^d \to \R$ such that the solution $u$ to 
				\begin{equation}
					\begin{cases}
						\Delta_{\R^d}^2 u(x) = f(x), \quad &x \in (-1,1)^d,\\
						u(x) = \pa_n u(x) = 0, \quad &x \in \pa (-1,1)^d.
					\end{cases}
				\end{equation}
				has $\| u \|_{C^5([-1,1]^d)} \leq C$. Let $\Lambda = \cl_1 \car_L \setminus \car_{L - \ell}$ be a boundary layer of width $\ell + 1$. Then we have
				\begin{equation} \label{e.olKbdry}
					\| {K}_L a \|^2_{L^2(\Lambda)} + \| a \|_{L^2(\Lambda)}^2 \leq \frac{C(\ell+1)}{L}(1 + \1_{d=2} \log L)
				\end{equation}
				for a constant $C$ dependent only on $d$ and $\| u \|_{C^5([-1,1]^d)}$.
			\end{prop}
			\begin{proof}
				We use $K_{L} a = a - \Delta G_{L} \ast f_L$ and treat these two terms separately. By the triangle inequality, we have
				$$
				\| \Delta G_{L} \ast f_L \|_{L^2(\Lambda)} \leq \| \Delta G_{L} \ast f_L - \Delta u_{L} \|_{L^2(\cl_1 \car_{L})} + \| \Delta u_{L} \|_{L^2(\Lambda)}.
				$$
				for $u_{L} = L^{-d/2 + 2} u(x/(L+2))$. We can bound the first term using \eref{bil_cont_approx} by $CL^{-1/2}$. Moreover, we have
				$$
				\Delta u_{L} = L^{-d/2} (\Delta_{\R^d} u)(x/(L+2)) + O(L^{-d/2 - 1}) \| u \|_{C^3([-1,1]^d)} \leq CL^{-d/2}
				$$
				and so $\| \Delta u_{L} \|^2_{L^2(\Lambda)} \leq C (\ell+1) L^{-1}$.
				
				We can bound $\| a \|_{L^2(\Lambda)}$ from the following estimate on $\| a \|_{L^\infty(\cl_1 \car_L)}$:
				$$
				|a(x)| \leq \sum_{y \in \car_L} \Gamma_L(x,y) f_L(y) \leq \| f_L \|_{L^\infty} \sum_{y \in \car_L} |\Gamma_L(x,y)| \leq CL^{-d/2}(1 + \1_{d=2} \log L).
				$$
				The last inequality follows from summing the inequality \eref{gam_sup} in \pref{greens_theo}. Putting all the above estimates together yields \eref{olKbdry}.
			\end{proof}
			
			We pause now to note that our condition on $f_L$ is certainly not sharp, and can be improved in dimensions $2$ and $3$ on square/cubic domains using the results of \cite{MS19}, one of which we now state. The following estimate on $G_L$ is readily seen to imply the result we need in the context of assumption \eref{A2} of \sref{concl_sec}. Actually, the paper \cite{MS19} proves very detailed estimates on $G_L$ far exceeding what we state and need here.
			
			\begin{theo}[\cite{MS19}, Theorem 1.1] \label{t.DeltaGthm}
				Let $d=2,3$. There is a constant $C$ such that
				$$
				|\Delta_x G_L(x,y)| \leq C(1 + |x-y|)^{2-d}(1 + \1_{d=2} \log L),
				$$
				uniformly in $L$.
			\end{theo}
			
			Consequently, for any $f_L \in \R^{\car_L}$, we have
			\begin{equation} \begin{aligned}
					\sup_{x \in \cl_1 \car_L} |\Delta (G_L \ast f_L)(x)| &\leq \sup_x \| \Delta G_L(x,\cdot) \|_{L^1(\cl_1 \car_L)} \cdot \sup_{x} |f_L(x)| \\ &\leq C L^2(1 + \1_{d=2} \log L) \sup_{x} |f_L(x)|.
			\end{aligned}\end{equation}
			In the case that $f_L = L^{-d/2 - 2} f(x/L)$ for some continuous $f$ bounded on $[-1,1]^d$ and ${a = -\Gamma_L \ast f_L}$, we have
			$$
			{K}_L^\perp a = \Delta G_L \Delta a = \Delta G_L f_L,
			$$
			and so
			\begin{equation}
				\sup_{x \in \cl_1 \car_L} |{K}_L^\perp a(x)| + |{K}_L a(x)| \leq C \| f \|_{L^\infty}(1 + \1_{d=2} \log L)L^{-d/2}.
			\end{equation} 
		\end{subsection}

	\end{section}

	%
	%
	
	\begin{acks}[Acknowledgments]
		The author was partially supported by NSF grant DMS-2000205. 
	\end{acks}
	\bibliographystyle{imsart-number} 
	\bibliography{bibliography}       
	
	
\end{document}